\documentclass[a4paper,12pt]{amsart}
\usepackage{amssymb}
\usepackage{amsmath, amsthm, amscd, amsfonts, amssymb, graphicx, color}
\usepackage{amsmath, amsthm, amscd, amsfonts, amssymb, graphicx, color}

\usepackage[cp1250]{inputenc}
\usepackage[T1]{fontenc}

\newtheorem{theorem}{Theorem}[section]
\newtheorem{lemma}[theorem]{Lemma}
\newtheorem{proposition}[theorem]{Proposition}
\newtheorem{corollary}[theorem]{Corollary}
\theoremstyle{definition}
\newtheorem{definition}[theorem]{Definition}
\newtheorem{example}[theorem]{Example}

\theoremstyle{remark}
\newtheorem{remark}[theorem]{Remark}
\numberwithin{equation}{section}

\begin{document}

	\vspace{15pt}
	
	\title{Disjoint $\mathcal{F}-$semi-transitivity in Banach algebras}
	\author{Stefan Ivkovi\' c}
	\address{Mathematical Institute of the Serbian Academy of Sciences and Arts, Kneza Mihaila 36, Beograd
		11000, Serbia}
	\email{stefan.iv10@outlook.com}
	\maketitle
	
		\begin{abstract}
		Motivated by the concept of disjoint topological transitivity for supercyclicity given in  \cite{novo-prvi}, in this paper, we consider the concept of disjoint Furstenberg-semi-transitivity for operators that are a composition of an isometric isomorphism and a left multiplier on a normed algebra. Thus, we characterize disjoint $\mathcal{F}$-semi-transitive and disjoint supercyclic such operators on a large class of non-unital normed algebras. It turns out that generalized weighted bilateral shifts on the standard Hilbert C*-module are just a special case of our theory. Generalized weighted composition operators on the normed algebra of operator-valued continuous functions vanishing at infinity on a locally compact, non-compact Hausdorff space are another special case of our theory. Next, we characterize disjoint $\mathcal{F}$-semi-transitive and disjoint supercyclic weighted composition operators on a large class of weighted solid Banach function spaces and we apply our results to the case of translations on weighted Morrey spaces. We illustrate all the results in this paper with concrete examples. 
	\end{abstract}
	
	
	\vspace{15pt}
	
	\begin{flushleft}
	\end{flushleft}
	
	\

	\section{Introduction}
	
	Supercyclicity is a fundamental concept in the dynamics of operators and numerous papers have been written on the topics related to supercyclicity of operators. For example,  Hilden and Wallen in \cite{Supercyclic-1} proved that any unilateral backward weighted shift is supercyclic. Afterwards, Salas in \cite{salas} characterized supercyclic bilateral weighted shift operators on $ l^p ( \mathbb{Z} ) $ in terms of a supercyclicity criterion. Supercyclicity of several kinds of operators has also been studied in for instance \cite{comparison, novo-prvi, novo-drugi, banach}. Moreover, quite recently, a characterization of supercyclic and topologically semi-transitive weighted composition operators on solid Banach function spaces, in particular Morrey spaces, and Segal algebras have been given in \cite{CAOT2}. Further, there is a close connection between supercyclicity and semi-Fredholm theory, see \cite{aiena, cao, Fredholm-3}. On the other hand, Furstenberg topological transitivity as a generalization of topological transitivity and hypercyclicity has been studied in several papers, see \cite{taba, furstenberg-prvi, furstenberg-drugi, furstenberg-treci, furstenberg-cetvrti}. Motivated by these facts, in this paper we consider the concept of $\mathcal{F}-$semi-transitivity (Definition \ref{glavna-definicija}), as a generalization of the concept of topological transitivity for supercyclicity given in \cite{novo-prvi} and we characterize disjoint $\mathcal{F}-$semi-transitive operators on a large class of non-unital normed algebra (Theorem \ref{algebra} and Corollary \ref{algebra-primena}). The operators we study are in fact a composition of an isometric isomorphism and a left multiplier on a non-unital normed algebra. The reason for considering these operators is because many known classes of operators that have earlier been studied in connection with the dynamics, such as bilateral weighted shifts on $l^2 ( \mathbb{Z} ) ,$ weighted composition operators on algebras of functions, generalized weighted shifts on Hilbert C*-modules and so on, turn out to be of this kind.
	
	Next, the dynamics of cosine operator functions has been studied in several papers, see \cite{ChangChen2013,Chen2014,Chen2015,Chen2016,Chen2023,tsi, FIL, Grada, Kalmes2010,Kostic2012, Wang2025}.   Although there are several published papers on the dynamics of cosine operator functions, supercyclicity of cosine operator functions has not been studied in any of the published papers, as far as we know. All these facts naturally lead us to study in this paper supercyclicity and $\mathcal{F}-$semi-transitivity (as a generalization of supercyclicity) of cosine operator functions generated by the above mentioned operators on normed non-unital algebras, i.e. the operators that are a composition of an isometric isomorphism and a left multiplier. In Proposition \ref{cosinus} and Corollary \ref{cosinus-gust} we provide some sufficient conditions for the cosine operators functions generated by these operators to be $\mathcal{F}-$semi-transitive on the respective normed algebra.
	
	Now, an attempt to axiomatically study hypercyclicity of operators on C*-algebras has been done in \cite[Section 4]{FIL}. This approach is limited to C*-algebras and does not cover many important and relevant cases in linear dynamics. The approach in Section 3 in this paper is more general because we consider a large class of non-unital normed algebras that are not necessarily C*-algebras. Our approach covers also the case of generalized bilateral weighted shift operators on standard Hilbert C*-modules defined in \cite{CAOT, FIL}.
	
	Section 4 of this paper is devoted to applications of the results from Section 3. Thanks to Corollary \ref{cosinus-gust}, in Proposition \ref{cosinus-shift} in this paper we provide sufficient conditions for $\mathcal{F}-$semi-transitivity of the cosine operator function generated by the generalized bilateral weighted shift operator on the standard Hilbert module over the C*-algebra of compact operators on a separable Hilbert space, extending thus \cite[Proposition 3.5]{FIL} from the case of topological transitivity to the case of $\mathcal{F}-$semi-transitivity. Example \ref{primer-modul-kompakt}  illustrates that the sufficient conditions of Proposition \ref{cosinus-shift} in this paper are strictly weaker than the sufficient conditions given in \cite[Proposition 3.5]{FIL}. Next, thanks to Theorem \ref{algebra} in this paper, in Corollary \ref{modul} we characterize disjoint $\mathcal{F}-$semi-transitive generalized bilateral weighted shift operators on the standard Hilbert module over a commutative non-unital C*-algebra, extending in this way \cite[Theorem 2.2]{CAOT} and \cite[Corollary 2.4]{CAOT} from the case of hypercyclicity and topological transitivity to the case of disjoint $\mathcal{F}-$semi-transitivity. Example \ref{primer-shift} shows that the conditions of Corollary \ref{modul} are strictly weaker than the conditions of \cite[Theorem 2.2]{CAOT} and \cite[Corollary 2.4]{CAOT}. Moreover, it is worth mentioning that Corollary \ref{modul} also covers the case of non-invertible generalized bilateral weighted shifts, whereas in \cite[Theorem 2.2]{CAOT} and \cite[Corollary 2.4]{CAOT}  invertibility of these shifts is assumed.  \\
	However, in addition of giving extensions of previously published results, in Section 4 we provide also a completely new application in the context of generalized weighted composition operators on the space of operator valued continuous functions. More precisely, we denote by $ \mathcal{C}$ the C*-algebra of compact operators on a separable Hilbert space $H ,$ and in Corollary \ref{operator-algebra} we characterize disjoint $\mathcal{F}-$semi-transitive generalized weighted composition operators on the algebra $C_{0}(\Omega,  \mathcal{C})$ equipped with the supremum-norm, i.e. the normed algebra of all $ \mathcal{C}-$valued continuous functions vanishing at infinity on a locally compact, non-compact Hausdorff space $ \Omega .$ The generalized weighted composition operators that we study in this context are in fact a composition of a weighted composition operator on  $C_{0}(\Omega,  \mathcal{C})$ and a right multiplication by a unitary operator on $H.$ Thus, these operators are at a same time a generalization of (ordinary) weighted composition operators on the classical space $ C_{0}(\Omega), $ and elementary (wedge) operators on $\mathcal{C}$ considered in \cite{BIMS}. We give also a concrete example where the conditions of Corollary \ref{operator-algebra} are satisfied (Example \ref{primer-operator-algebra}).
	
	As noticed in Remark \ref{Banach-bimodule} in this paper, the results in Section 3 can under certain conditions be extended from normed algebras to Banach bimodules over normed algebras. We study this case more in detail in Section 5. Motivated by the concept of solid Banach function spaces from \cite{saw} as a generalization of Lebesgue, Orlicz and Morrey spaces, in Section 5 we study the dynamics of weighted composition operators on these spaces. By considering these spaces as bimodules over algebras of bounded measurable functions and by applying Theorem \ref{algebra} together with Remark \ref{Banach-bimodule}, we characterize  disjoint $\mathcal{F}-$semi-transitive weighted composition operators on these solid spaces. However, we also go a step further and in Theorem \ref{vekt}  we characterize  disjoint $\mathcal{F}-$semi-transitive weighted composition operators on weighted solid Banach function spaces. The motivation for considering weighted solid spaces comes from \cite{banach} where disjoint supercyclicity of weighted translations on weighted Orlicz spaces has been studied. Since weighted solid Banach function spaces are a generalization of weighted Orlicz spaces, it follows that Theorem \ref{vekt} in this paper  is a generalization of the results in \cite{banach}. Theorem \ref{vekt} covers also the case of weighted Morrey spaces, so in Example \ref{zavrsetak} we apply Theorem \ref{vekt} in the concrete case of translations on weighted Morrey spaces. Further, Theorem \ref{vekt} is an extension of \cite[Proposition 3.5]{CAOT2} from supercyclicity on solid Banach function spaces to disjoint $\mathcal{F}-$semi-transitivity on weighted solid spaces. In particular, from \cite[Proposition 3.5]{CAOT2} it follows that non-weighted translations on Morrey spaces can never be supercyclic, however, by Theorem \ref{vekt} and Example \ref{zavrsetak} we conclude that there exist weighted Morrey spaces where non-weighted translations are supercyclic and $\mathcal{F}-$semi-transitive, showing the difference between weighted and non-weighted solid Banach function spaces.
	
	Our opinion is that the results from Section 3 in this paper has a potential to be applied also to some other concrete examples of normed algebras and Banach bimodules, so in Section 6 we leave this as an open question for further research.
	\section{Preliminaries}
	We need first the following definition.
	\begin{definition}\label{glavna-definicija}
		Let $\mathcal{S}$ be a set, $\mathcal{F}$ be a family of subsets of $\mathcal{S} $   
		and $$\{ T_{t,1} \}_{t \in \mathcal{S}}, \dots,
		\{ T_{t,N} \}_{t \in \mathcal{S}}$$ be families of bounded linear operators on a Banach space $X$.  
		We say that $\{ T_{t,1} \}_{t \in \mathcal{S}}, \dots, \{ T_{t,N} \}_{t \in \mathcal{S}}$ are disjoint $\mathcal{F}$-semi-transitive, or shortly $d\mathcal{F}$-semi-transitive, if  
		for every collection of non-empty open subsets $\mathcal{O}, V_{1}, \dots, V_{N}$ of $X$,  
		there exists some $ F \in \mathcal{F} $  such that for all  $ t \in F$ there exists some $  \lambda_{t} \in \mathbb{R}^+  $  satisfying that 
		$$
		\mathcal{O} \cap \lambda_{t} T_{t,1}^{-1} (V_{1}) \cap \dots \cap \lambda_{t} T_{t,N}^{-1}(V_{N}) \neq \varnothing.
		$$
		
	\end{definition}
	We recall also that the families  $\{ T_{t,1} \}_{t \in S}, \dots, \{ T_{t,N} \}_{t \in S}$ are said to be \textit{disjoint $\mathcal{F}$-transitive, or shortly $d\mathcal{F}$-transitive}, if  
	for every collection of non-empty open subsets $\mathcal{O}, V_{1}, \dots, V_{N}$ of $X$,  
	there exists some $ F \in \mathcal{F} $  such that for all  $ t \in F$ we have
	$$
	\mathcal{O} \cap  T_{t,1}^{-1} (V_{1}) \cap \dots \cap  T_{t,N}^{-1}(V_{N}) \neq \varnothing.
	$$
	
	\text{ }
	
	Let now $\mathcal{A}$ be a non-unital normed algebra such that $\mathcal{A}$ is a left ideal in a unital normed algebra $\mathcal{A}_1$.
	We do \underline{not} assume that the norm on $\mathcal{A}$ extends to $\mathcal{A}_1$.
	We will say that $\mathcal{A}$ satisfies \textit{ the condition (E) } with respect to $\mathcal{A}_1$ if $$ \| ba \| \leq \| b \|_1 \| a \| $$ for all $ a \in \mathcal{A}$ and $ b \in \mathcal{A}_1 ,$ where $ \| \text{ } \|_1 $ denotes the norm on $\mathcal{A}_1 .$\\
	Let  $\{p_\alpha\}_\alpha$ be a set in $\mathcal{A}$ satisfying that given any open subset $O$ of $\mathcal{A}$ and $ x \in O$ there exists some $ p_{\alpha_0} \in \{p_\alpha\}_\alpha$ such that $ p_{\alpha_0}^{3} x \in O$  as well. This condition will be called \textit{the condition} $ (\mathcal{P}.) $ \\

	For a set $\mathcal{S},$ let $\mathcal{F}$ be a family of subsets of $\mathcal{S}$, and let
	$$
	\{\Phi_{t,1}\}_{t\in \mathcal{S}},\ \{\Phi_{t,2}\}_{t\in \mathcal{S}},\ \ldots,\ \{\Phi_{t,N}\}_{t\in \mathcal{S}}
	$$
	be families of isometric algebra isomorphism of $\mathcal{A}_1$ such that
	$$
	\Phi_{t,j}(\mathcal{A}) = \mathcal{A}
	\qquad (t\in \mathcal{S},\ 1\le j \le N).
	$$
	and such that $\Phi_{t,j\vert \mathcal{A}}$ is an isometry for all $j \in \lbrace 1,\dots , N \rbrace.$ We will assume that the system
	$$
	\{\Phi_{t,1}\}_{t\in \mathcal{S}},\ \{\Phi_{t,2}\}_{t\in \mathcal{S}},\ \ldots,\ \{\Phi_{t,N}\}_{t\in \mathcal{S}}
	$$
	is disjoint aperiodic with respect to $\{p_\alpha\}_\alpha$, that is for each fixed $p_{\alpha}$ and every $H\in \mathcal{F}$ there exist some $F \subseteq H$ with $F \in \mathcal{F}$ such that
	$$
	p_{\alpha} \, \Phi_{t,\ell}(p_{\alpha}) = 0, 
	\quad 
	p_{\alpha} \, \Phi_{t,\ell}^{-1}(p_{\alpha}) = 0, 
	\quad \text{and} \quad
	p_{\alpha} \, \Phi_{t,\ell}\!\big( \Phi_{t,r}^{-1}(p_{\alpha}) \big) = 0
	\quad \text{for all }  t \in F
	$$
	$ \text{ and all } r, \ell \in \{1, \dots, N\} \text{ with } r \neq \ell.
	$
	Further, we will assume that for each $\alpha$ and all $t \in \mathcal{S}$, 
	$r, \ell \in \{1, \dots, N\}$, and $a \in \mathcal{A}$, 
	it holds that
	$$
	\| a \, p_{\alpha} \| \le \| a \|, 
	\quad 
	\| a \, \Phi_{t,\ell}(p_{\alpha}) \| \le \| a \|,
	$$
	$$
	\| a \, \Phi_{t,\ell}^{-1}(p_{\alpha}) \| \le \| a \|,
	\quad 
	\quad 
	\| a \, \Phi_{t,\ell}\big(\Phi_{t,r}^{-1}(p_{\alpha})\big) \| \le \| a \|. 
	$$
	This condition will be called \textit{ the condition (R) } throughout this section.
	
	We now set up families 
	$\{ b_{t,1} \}_{t \in \mathcal{S}}, \cdots,\{ b_{t,N} \}_{t \in \mathcal{S}}$
	of elements in $\mathcal{A}_1$, 
	and for each $t \in \mathcal{S}$ and $\ell \in \{1, \dots, N\}$, 
	we define the operator
	$$
	T_{t,\ell} : \mathcal{A} \to \mathcal{A}
	\quad \text{by} \quad
	T_{t,\ell}(a) = b_{t,\ell} \, \Phi_{t,\ell}(a).
	$$
	
	Due to the condition $(E)$, 
	it is not hard to check that $T_{t,\ell}$ 
	is a bounded linear operator on $\mathcal{A}$.  
	Further, we will assume that for each $ \alpha $ there exists so called ''$ \alpha$-inverse'' for $b_{t,\ell} $ denoted by $ b_{t,\ell, \alpha}^{-1} ,$  that is for each $ \alpha $ there exists an element $ b_{t,\ell, \alpha}^{-1}  \in \mathcal{A}_1$ satisfying that $ b_{t,\ell} b_{t,\ell, \alpha}^{-1} a p_\alpha = b_{t,\ell, \alpha}^{-1} b_{t,\ell} a p_\alpha = a p_\alpha$ for all $ a \in \mathcal{A} .$
	
	 Also, we will assume that for each $ t \in \mathcal{S} $ and $\ell \in \{1, \dots, N\}$ it holds that   $ b_{t,\ell} a = 0 $ if and only if $ a = 0.$ Such condition will be called \textit{the condition (C)} throughout the manuscript.
	 
	\section{Main results}
	
	We have the following theorem.
	
	\begin{theorem}\label{algebra}
		Under the above notation and assumptions, the following statements are equivalent.
		\begin{enumerate}
			\item The families $\{ T_{t,1} \}_{t \in \mathcal{S}}, \cdots,\{ T_{t,N} \}_{t \in \mathcal{S}}$ are $d \mathcal{F}$-semi-transitive.
			
			\item For each fixed $ p_\alpha $ and $ \varepsilon > 0 $ there exists some 
			$ F \in  \mathcal{F} $ and families $\{ d_t \}_{t \in F}$, $ \{ g_{t,1} \}_{t \in F} \dots , \{ g_{t,N} \}_{t \in F} $ in $ \mathcal{A} $
			such that 
			$
			\|d_t - p_\alpha^{2}   \| < \varepsilon, \quad    \|g_{t,l} - p_\alpha^2   \| < \varepsilon,
			$
			and 
			$$
			\| b_{t,r, \alpha}^{-1} g_{t,r}    \|     \| \Phi_{t,l}^{-1} (b_{t,l}) d_{t}   \| < \varepsilon^{2}
			\quad \text{for all } t \in F \text{ and } r,l \in \{ 1, \ldots, N \}.
			$$
			Moreover, for each distinct $ r,l \in \{ 1, \ldots, N \}$ and $t \in F$ it holds that
			$$
			\| \Phi_{t,r}^{-1}(b_{t,r})   \Phi_{t,l}^{-1}(b_{t,l, \alpha}^{-1}  g_{t,l})   \| < \varepsilon.
			$$
		\end{enumerate}
		
	\end{theorem}

	\begin{proof}
		We prove first \textit{(1) $\Rightarrow$ (2)}. Since the families $\{ T_{ \Phi_{t,1},b_{t,1} } \}_{t \in S}  , \dots ,   \{ T_{ \Phi_{t,N},b_{t,N} } \}_{t \in S}$ 
		are $d \mathcal{F}$-semi-transitive, given $\varepsilon > 0$ and some fixed $ p_{\alpha} $ we can find 
		$a_t \in \mathcal{A}$ and $\lambda_{t}  \in \mathbb{R}^{+}$  such that $   \| a_t - p_{\alpha}    \| < \varepsilon$ and 
		$
		\|\lambda_{t} T_{ \Phi_{t,l}},b_{t,l}  (a_t)  -p_{\alpha}       \| < \varepsilon   \text{ for all } t \in H \text{ and some } H \in \mathcal{F}.
		$
		This gives
		$
		\bigl   \|\,(a_{t}-p_{\alpha})\,p_{\alpha}\,\bigr   \|\;\le\;\bigl   \|a_{t}-p_{\alpha} \bigr   \|\;<\;\varepsilon.
		$

		Moreover, since
		$
		\bigl\{\    {\Phi_{t,l}} \bigr\} _{  t \in S, \; 1\le l \le N }          
		$
		is a disjoint aperiodic system, we can find some $ F \in \mathcal{F}$ with $ F \subseteq H $ such that 
		$
		p_{\alpha}  \!   \Phi_{t,l}^{-1}( p_{\alpha} )  =0, \;  
		p_{\alpha}  \!   \Phi_{t,l}( p_{\alpha} )  =0, \;  
		$
		\noindent and
		$
		\Phi_{t,r}( \Phi_{t,l}^{-1}( p_{\alpha} ) ) =0, \; \text{for all } \; t\in F \; \text{ and } r,l \in \{ 1, \dots N\} \text{ with } r \neq l.
		$
		Hence, we get 
		
		$$
		\qquad
		\frac{1}{|\lambda_t|}\,
		\left   \|
		\Phi_{t,r}^{-1}\!\bigl(b_{t,r, \alpha}^{-1}\bigr)\,
		\Phi_{t,r}^{-1}\!\bigl(b_{t,r}\bigr)\, a_{t}\, \Phi_{t,r}^{-1}(p_{\alpha})\, \lambda_{t}
		\right   \|
		= \left   \|
		\! b_{t,r, \alpha}^{-1}\,
		\!\bigl(b_{t,r}\, \Phi_{t,r} (a_{t})\, (p_{\alpha})\, \bigr) \lambda_{t} 
		\right   \|
		$$
		$$
		= \left   \|
		a_t\,\Phi_{t,r}^{-1}(p_{\alpha})
		\right   \| =\left   \|(a_{t}-p_{\alpha})\,\Phi_{t,r}^{-1}(p_{\alpha})
		\right   \|
		\le 
		\left   \|a_{t}-p_{\alpha}\right   \|
		< \varepsilon .
		$$ and 
		$$
		\left   \|
		\lambda_t\, b_{t,r}\,\Phi_{t,r}(a_{t})\,p_{\alpha}
		\;-\; p_{\alpha}^{2}
		\right   \|
		\le
		\left   \|
		\lambda_t\, b_{t,r}\,\Phi_{t,r}(a_t) - p_{\alpha}
		\right   \|
		=
		\left   \|
		\lambda_{t}\,T_{t,r}(a_{t}) - p_{\alpha}
		\right   \|
		< \epsilon,
		\qquad 
		$$ for all  $t\in F$ and $ r\in\{1,\dots,N\}$. \\
		For each $t\in F$ and $r\in\{1,\dots,N\},$  we put $	g_{t,r} = \lambda_t\, T_{t,r}(a_t) p_\alpha.$ 
		Then $$
		\|g_{t,r}-p_{\alpha}^{2}   \| \;<\; \varepsilon$$ and
		$$
		\frac{1}{|\lambda_t|}\,
		\bigl   \|\Phi_{t,r}^{-1}\!\bigl(b_{t,r, \alpha}^{-1}\bigr)\,\Phi_{t,r}^{-1}(g_{t,r})\bigr)\bigr   \|
		\;=\; \frac{1}{|\lambda_t|}\,   \|b_{t,r, \alpha}^{-1} g_{t,r}   \|
		\;<\; \varepsilon .$$
		
		Similarly, for each $ t\in F $ and	$ l\in\{1,\dots,N\},$ we have 
		$$
		\|\lambda_t\, b_{t,l}\,\Phi_{t,l}(a_t)\,\Phi_{t,l}(p_\alpha)   \|
			\;=\;    \|(\lambda_t\, b_{t,l}\,\Phi_{t,l}(a_t)-p_\alpha)\,\Phi_{t,l}(p_\alpha)   \|
			\;\le\;    \|\lambda_t\, T_{t,l}(a_t)-p_\alpha   \|
			\;<\; \varepsilon ,$$
		
		and, as shown above,	
		$$ \| a_t p_\alpha - p_{\alpha}^{2}   \| \;<\; \varepsilon $$
			For each $ t\in F$ and $ l\in\{1,\dots,N\},$
			 set $ d_t = a_t p_\alpha. $
			Then   $$  \|d_t-p_\alpha^{2}   \| \;<\; \varepsilon $$
			and 
			$$
			|\lambda_t|\,\bigl   \|\Phi_{t,l}^{-1}(b_{t,l})\,d_t\bigr   \|
			\;=\;    \|\lambda_t\, b_{t,l}\,\Phi_{t,l}(a_t p_\alpha)   \|
			\;<\; \varepsilon .
		$$
		% --- end snippet ---

		Hence, for all $t\in F$ and $ r,l\in\{1,\dots,N\},$  we have
		$$\| b_{t,r, \alpha}^{-1} (g_{t,r})    \|     \| \Phi_{t,l}^{-1} (b_{t,l}) d_{t}   \| < \varepsilon^{2}. $$

		Further, for all $t\in F$ and each distinct $r,l\in\{1,\dots,N\}$  we get
		$$
		\Bigl\|\,\Phi_{t,r}^{-1}\!\bigl(b_{t,r}\bigr)\, \Phi_{t,l}^{-1}\!\bigl(b_{t,l, \alpha}^{-1}\bigr)
		\Phi_{t,l}^{-1}\!\bigl(g_{t,l}\bigr)\Bigr\|
		\;=\; \Bigl\|\,\!b_{t,r}\, \Phi_{t,r} \Phi_{t,l}^{-1}\!\bigl(b_{t,l, \alpha}^{-1}
		g_{t,l}\bigr)\Bigr\|
		$$
		$$
		=\bigl\|\,\lambda_t\, b_{t,r}\,\Phi_{t,r}(a_t)\,
		\Phi_{t,r}\!\bigl(\Phi_{t,l}^{-1}(p_\alpha)\bigr)\bigr\| 
		$$
		$$
		=\; \bigl\|\,(\lambda_t\, T_{t,r}(a_t)-p_\alpha)\,
		\Phi_{t,r}\!\bigl(\Phi_{t,l}^{-1}(p_\alpha)\bigr)\bigr\|
		\;<\; \varepsilon ,
		$$
		
		\noindent
		where we have used that
		$
		p_\alpha\,\Phi_{t,r}\!\bigl(\Phi_{t,l}^{-1}(p_\alpha)\bigr)=0
		\quad\text{for all } t\in F
		\text{ and each distinct } r,l\in\{1,\dots,N\}.
		$

		Now we prove the implication $(2)\Rightarrow(1)$. 
		Let $O,V_1,\dots,V_N$ be non-empty open subsets of $\mathcal{A}$. 
		Then we can find some non-zero $x\in O$ and non-zero $y_l\in V_l$ for each 
		$l\in\{1,\dots,N\}$. 
		There exists some  $\alpha$ such that
		$ p_\alpha^{3} x\in O $
		and
		$ p_\alpha^{3}y_l\in V_l\setminus\{0\} $
		for each $l=1,\dots,N$. Clearly, this gives that $ p_\alpha x \neq 0 $ and $ p_\alpha y_l \neq 0 $ for each $l=1,\dots,N$. Also, since $O\setminus\{0\},\,V_1\setminus\{0\},\dots, V_N\setminus\{0\}$ are open and non-empty,
		there is some $\delta>0$ such that the $\delta$-neighbourhood of $p_\alpha^{3}x$ is
		contained in $O \setminus\{0\}$, and the $\delta$-neighbourhood of $p_\alpha^{3}y_l$ is contained
		in $V_l\setminus\{0\}$ for each $l=1,\dots,N$.
		
		Let
		$ C_y=\max_{1\le l\le N}\bigl\{\|y_l\|,\bigr\} $
		and set
		$$ \varepsilon
		:=\min\!\left\{
		\frac{\delta}{4\|x\| \| p_\alpha    \|},\;
		\frac{\delta}{4N\sqrt{\|x\|C_{y}} \| p_\alpha    \| } 
		\frac{\delta}{4 N C_{y}   \| p_\alpha    \|} , \frac{\| p_{\alpha}^{3} x   \| }{2 \| p_{\alpha} x   \|},\frac{\| p_{\alpha}^{3} y_1   \| }{2 \| p_{\alpha} y_1   \|}, \cdots, \frac{\| p_{\alpha}^{3} y_N   \| }{2 \| p_{\alpha} y_N   \|}
		\right\}. $$
		
		Choose $F\in  \mathcal{F}$ and the families 
		$ \{ g_{t,1} \}_{t\in F}, \dots  , \{ g_{t,N} \}_{t\in F}, \{d_{t} \}_{t \in F} $ 
		satisfying the assumptions in $(2)$ with respect to $ \epsilon $ and $p_{\alpha} $.For each $ \alpha $ and every $ t \in F, l\in\{1,\dots,N\}  $ we define the operator $ S_{t,l, \alpha}$ on $ \mathcal{A} $ by $ S_{t,l, \alpha}(a) = \Phi_{t,l}^{-1}(b_{t,l, \alpha}^{-1}) \Phi_{t,l}^{-1}(a) $ for all $  a \in \mathcal{A} .$ We observe now that for each $t\in F$ and $l\in\{1,\dots,N\}$ it holds that
		$$
		\parallel T_{t,l}\!\bigl(d_{t} p_\alpha x\bigr)\parallel
		= \parallel\, b_{t,l}\,\Phi_{t,l} \bigl(d_{t}\bigr)\,\Phi_{t,l}( p_\alpha x)\,\parallel
		\le \parallel\,  b_{t,l} \, \Phi_{t,l} \bigl(d_{t}\bigr)\,\parallel\,
		\parallel\,\Phi_{t,l}(p_\alpha x)\,\parallel
		$$
		$$
		= \parallel\, \Phi_{t,l}^{-1} \, (b_{t,l}) \, d_{t}\,\parallel\;
		\parallel \, x \,\parallel \| p_\alpha    \|
		$$
		and
		$$
		\parallel S_{t,l, \alpha}\!\bigl(g_{t,l} p_\alpha \,y_{ \ell}\bigr)\parallel
		= \parallel\, \Phi_{t,l}^{-1}\, (b_{t,l, \alpha}^{-1} )\,\Phi_{t,l}^{-1}\!\bigl(g_{t,l}\bigr)\,
		\Phi_{t,l}^{-1}(p_\alpha \,y_{ \ell})\,\parallel
		$$
		$$
		\le \parallel\, \Phi_{t,l}^{-1}\,(b_{t,l, \alpha}^{-1} \, g_{t,l} )  \,\parallel\,
		\parallel\, \Phi_{t,l}^{-1}(p_\alpha \,y_{ \ell})\,\parallel
		= \parallel\, b_{t,l, \alpha}^{-1}\, g_{t,l}\,\parallel\;
		\parallel\, y_{\ell}\,\parallel \| p_\alpha    \| .
		$$
		
		Moreover, for all $t\in F$ and each distinct $r,l \in \lbrace 1,\dots,N \rbrace$ we have
		$$
		\parallel T_{t,r}\!\bigl(S_{t,l, \alpha}(g_{t,l}p_\alpha \,y_{ \ell})\bigr)\parallel
		= \parallel\,  b_{t,r}\, \Phi_{t,r}\, ( \Phi_{t,l}^{-1}\, (b_{t,l, \alpha}^{-1}) \bigr) \,
		\Phi_{t,r}\, \bigl( \Phi_{t,l}^{-1}\, (g_{t,l} p_\alpha \,y_{ \ell}\bigr) \bigr)\,\parallel
		$$
		$$
		= \parallel 
		\Phi_{t,r}^{-1}\, ( b_{t,r}\,)   \Phi_{t,l}^{-1}\, (b_{t,l, \alpha}^{-1}) \,  \Phi_{t,l}^{-1}\, (g_{t,l}) \, \Phi_{t,l}^{-1}\, (p_\alpha \,y_{ \ell})\,\parallel
		$$
		$$
		\leq \parallel 
		\Phi_{t,r}^{-1}\, ( b_{t,r}\,)   \Phi_{t,l}^{-1}\, (b_{t,l, \alpha}^{-1}g_{t,l}) \, \,\parallel \, \parallel y_{\ell} \parallel \| p_\alpha    \| .
		$$
		Finally, since $ b_{t,l, \alpha}^{-1} $ is the $ \alpha$-inverse of $b_{t,\ell} ,$ it follows that $$ T_{t,l}\!\bigl(S_{t,l, \alpha}(g_{t,l}p_\alpha \,y_{ \ell})) = g_{t,l}p_\alpha \,y_{ \ell} $$ for all $t\in F$ and each $ \ell \in \lbrace 1,\dots,N \rbrace .$

		Since $\| d_t p_{\alpha} x -  p_{\alpha}^{3} \, x \| < \| \epsilon \| \| p_{\alpha} x   \| < \|  p_{\alpha}^{3}  \, x \|$, 
		it follows that $d_t p_{\alpha} x \neq 0.$ Similarly, $g_{t,l}p_{\alpha} y_{\ell} \neq 0$  for all
		$t \in F$ and $\ell \in \{1, \dots, N\}$. 
		Hence, since by the assumption we have that for each $ t \in \mathcal{S} $ and $\ell \in \{1, \dots, N\}$ it holds that   $ b_{t,\ell} a = 0 $ if and only if $ a = 0,$ it follows that $T_{t, \ell}(d_t p_\alpha x) \neq 0 $ for each $ t \in F$ and $\ell \in \{1, \dots, N\}.$  Similarly, since $ b_{t,l, \alpha}^{-1} $ is the $ \alpha$-inverse of $b_{t,\ell} ,$  it folllows that   $S_{t, \ell, \alpha}(g_{t,\ell} p_\alpha y_\ell) \neq 0$ for all $t \in F$ and $\ell \in {1, \dots, N}.$

		For each $t\in F$, we put
		$$
		x_{t} := d_{t} p_\alpha x + \frac{\sqrt{\sum_{l=1}^{N}\,\parallel T_{t,l}(d_t p_\alpha x)\parallel}}{\sqrt{\sum_{l=1}^{N}\,\parallel S_{t,l, \alpha}(g_{t,l} p_\alpha y_{\ell})\parallel} \qquad } \, {\sum_{l=1}^{N}\, S_{t,r, \alpha}(g_{t,r} p_\alpha y_{r})}.
		$$
		
		%	\medskip
		By triangle inequality and choice of $ \varepsilon$ one can deduce that for each $t\in F$, we have
		$$
		x_{t} \in O, \ \text{and} \ 
		\frac{
			\sqrt{ \sum_{\ell = 1}^{N} \| S_{t, \ell, \alpha}( g_{t, \ell} p_\alpha y_\ell ) \| }
		}{
			\sqrt{ \sum_{\ell = 1}^{N} \| T_{t, \ell}( d_t p_\alpha x ) \| }
		}
		T_{t, r}(x_t) \in V_r.
		$$
		for all $r \in \lbrace 1, \dots , N \rbrace.$ 
	\end{proof}

	The next corollary can easily be deduced from Theorem \ref{algebra}. 
	
	\begin{corollary}\label{algebra-primena}
		Let now $\mathcal{F}$ be the family of all infinite subsets of $\mathbb{N} .$ If there exist dense subsets $ \mathcal{A}_0, \mathcal{B}_1, \cdots, \mathcal{B}_N$ of $ \mathcal{A} $ and a strictly increasing sequence $\{ n_{k} \}_{k} \subseteq \mathbb{N}$ such that  
		for each $s, \ell \in \{1, \dots, N\},$ every fixed $ \alpha $  and all $ x \in \mathcal{A}_0, y_s \in \mathcal{B}_s $ we have $$ \lim_{k \to \infty} \| b_{n_k,s, \alpha}^{-1} y_s    \|     \| \Phi_{n_k,l}^{-1} (b_{n_k,l}) x   \| = 0,$$ and, in addition, for each distinct $ s, \ell \in \{ 1, \ldots, N \}$ it holds that
		$$
		\lim_{k \to \infty}	\| \Phi_{n_k,\ell}^{-1}(b_{n_k,\ell})   \Phi_{n_k,s}^{-1}(b_{n_k,s, \alpha}^{-1}  y_{s})  \| =0,
		$$ then the families $\{ T_{n,1} \}_{n \in \mathbb{N}}, \cdots,\{ T_{n,N} \}_{n \in \mathbb{N}}$ are $d \mathcal{F}$-semi-transitive.
	\end{corollary}
	\begin{remark}
			We notice that the assumption that $\{ \Phi_{t, 1} \}_{t \in \mathcal{S}}, \cdots, \{ \Phi_{t, N} \}_{t \in \mathcal{S}}$ 
		is a disjoint aperiodic system in Theorem \ref{algebra} is only needed for the proof of the implication $(1) \Rightarrow (2)$. Also, we notice that we have used the condition \textit{(R)} only in the proof of the implication $(1) \Rightarrow (2)$ to obtain the desired inequalities. Hence, the assumption that $\{ \Phi_{t, 1} \}_{t \in \mathcal{S}}, \cdots, \{ \Phi_{t, N} \}_{t \in \mathcal{S}}$ 
		is disjoint aperiodic and the condition \textit{(R)} are not needed in Corollary \ref{algebra-primena}.
	\end{remark}
	
	Now we set $N = 1$ and we will therefore write $b_{t}$, $\Phi_{t}$ and $T_{t}$  instead of 
	$b_{t,1}$, $\Phi_{t,1}$ and $T_{t,1},$ respectively. We will now assume that $b_{t}$ is invertible 
	for each $t \in \mathcal{S}$ and we will consider the map 
	$S_{t} : \mathcal{A} \rightarrow  \mathcal{A} $ given by $S_{t}(a) = \Phi_{t}^{-1}(b_{t}^{-1}) \, \Phi_{t}^{-1}(a)$ for all 
	$a \in  \mathcal{A} $. Obviously, $ S_{t} = T_{t}^{-1} $ for each $t \in \mathcal{S}.$

	We set $C_{t} := \frac{1}{2}(T_{t} + S_{t})$ for all $t \in \mathcal{S}$.
	
	\begin{proposition}\label{cosinus}
		Under the above notation and assumptions, we have $(2) \Rightarrow (1)$.
		
		(1) The family $\{ C_{t} \}_{t \in \mathcal{S}}$ is $\mathcal{F}$–semi–transitive.
		
		(2) For each fixed $p_{\alpha}$ and $\varepsilon > 0$ there exists some $F \in \mathcal{F}$ 
		and families $\{ d_{t} \}_{t \in F}$, $\{ g_{t} \}_{t \in F}$, $\{ e_{t} \}_{t \in F}$, $\{ f_{t} \}_{t \in F}$ 
		with $e_{t} + f_{t} = g_{t}$ for all $t \in F$ such that 
		$$\| d_{t} - p_{\alpha}^{2}  \| < \varepsilon, \qquad 
		\| g_{t} - p_{\alpha}^{2} \| < \varepsilon,$$ 
		$$\| \Phi_{t}^{-2}  (b_{t})  \Phi_{t}^{-1}(b_{t}) e_{t} \| < \varepsilon,$$
		$$\|  \Phi_{t}(b_{t}^{-1}) \, b_{t}^{-1} \, f_{t} \| < \varepsilon,$$
		$$\|  \Phi_{t}^{-1}(b_{t}) \, d_{t} \| \, \|  \Phi_{t}^{-1} (b_{t}) \, e_{t} \|      < \varepsilon^{2},$$
		$$\|  \Phi_{t}^{-1}(b_{t}) \, d_{t} \| \, \| b_{t}^{-1} f_{t} \| < \varepsilon^{2},$$
		$$\|  b_{t}^{-1} \, d_{t} \| \, \| \Phi_{t}^{-1}(b_{t}) \, e_{t} \| < \varepsilon^{2}, $$
		$$\|  b_{t}^{-1} \, d_{t} \| \, \| b_{t}^{-1} \, f_{t} \| < \varepsilon^{2},  \text{ for all } t \in F.$$
	\end{proposition}

	\begin{proof}

		Given non–empty open subsets $\mathcal{O}$ and $\mathcal{V}$ of $\mathcal{A}$, as in the proof of Theorem \ref{algebra} part $(2)\Rightarrow(1)$
		we can find some $\delta > 0$ and some $p_{\alpha}^{3} x \in \mathcal{O} \setminus \lbrace 0 \rbrace$,
		$p_{\alpha}^{3} y \in \mathcal{V} \setminus \lbrace 0 \rbrace $ such that $\delta$–neighbourhoods of
		$p_{\alpha}^{3} x $ and $p_{\alpha}^{3} y $ are contained in $\mathcal{O} \setminus \{0\}$
		and $V \setminus \{0\}$, respectively.
		
		Set 
		$$\varepsilon := \min \left\{ \frac{\delta}{4 \| p_{\alpha} \| (\| x \| + \| y \|)}, \;
		\frac{\delta}{16 \| p_{\alpha} \| \sqrt{(\| x \| \| y \|)}}, \;
		\frac{\| p_{\alpha}^{3} x \|}{2 \| p_{\alpha} x \|}, \;
		\frac{\| p_{\alpha}^{3} y \|}{2 \|p_{\alpha} y \|} \right\} .$$
		
		Choose $F \in \mathcal{F}$ and families $\{ d_{t} \}_{t \in F}$, $\{ g_{t} \}_{t \in F}$ ,  $\{e_{t} \}_{t \in F}$, $\{ f_{t} \}_{t \in F}$ 
		satisfying the assumptions in (2) with respect to $p_{\alpha}$ and $\varepsilon$. Since $e_{t} + f_{t}=g_{t}$ and $g_{t} \,  p_{\alpha} y \neq 0$ 
		(as shown in the proof of  Theorem \ref{algebra} part $(2)\Rightarrow(1)$), it 
		follows that either $e_{t} \,  p_{\alpha} y \neq 0$ or $f_{t} \, p_{\alpha}  y \ne 0$. 
		Then, by the same arguments as in the 
		proof of Theorem \ref{algebra} part $(2)\Rightarrow(1)$ we can show that either 
		$T_{t}(e_{t}, p_{\alpha} y) \ne 0$ 
		or 
		$S_{t}(f_{t}, p_{\alpha} y) \ne 0$.\\
		Further, we notice that for each $t \in F$ it holds that
		$$\| T_{t}^{2}(e_{t}, p_{\alpha}  y) \| = 
		\| b_{t} \, \Phi_{t} \, ( b_{t} ) \, \Phi_{t}^{2} \, ( e_{t} \, p_{\alpha} y) \| 
		$$
		$$
		\le 
		\| \Phi_{t}^{-2} \, (b_{t}) \, \Phi_{t}^{-1} \, (b_{t})  \, e_{t} \| \, 
		\|  p_{\alpha} \|  \, \| y \|$$
		
		and,
		
		$$
		\| S_{t}^{2}(f_{t}, p_{\alpha} y) \| 
		= 
		\| \Phi_{t}^{-1}(b_{t}^{-1}) \cdot \Phi_{t}^{-2}(b_{t}^{-1} f_{t} \,  p_{\alpha} y ) \|
		\le 
		\| \Phi_{t} (b_{t}^{-1}) b_{t}^{-1} f_{t} \| 
		\, \|p_{\alpha} \| \, \|  \| y \|.
		$$
		
		For each $t \in F,$ put 
		$$v_{t} := d_{t} p_{\alpha}x + 
		\frac{2\sqrt{\| T_{t}(d_{t}  p_{\alpha} x) \| + \| S_{t}(d_{t} p_{\alpha} x) \|}}{\sqrt{\| T_{t}(e_{t} p_{\alpha} y) \| + \| S_{t}(f_{t} p_{\alpha} y) \|}}( T_{t}(e_{t} + p_{\alpha} x)  +  S_{t}(e_{t} p_{\alpha} x)  ) .
		$$
		By the above inequalities and inequalities in the proof of Theorem \ref{algebra} part $(2)\Rightarrow(1)$, 
		the choice of $\varepsilon$ and the triangle inequality, 
		we can deduce that $v_{t} \in \mathcal{O}$ and
		
		$$
		\frac{\sqrt{\| T_{t}(e_{t}  p_{\alpha} y) \| + \| S_{t}(f_{t} p_{\alpha} y) \|}}{\sqrt{\| T_{t}(d_{t} p_{\alpha} x) \| + \| S_{t}(d_{t} p_{\alpha} x) \|}}
		C_{t}(v_{t}) \in \mathcal{V} \text{ for all } t \in F.
		$$
	\end{proof}

	By putting either $e_{t} = 0$ (so that $f_{t} = g_{t}$), or by putting 
	$f_{t} = 0$ (so that $e_{t} = g_{t}$) in Proposition \ref{cosinus}, 
	then we can get the following corollary.
	
	\begin{corollary}\label{cosinus-gust}
		Let now $\mathcal{F}$ be the family of all 
		infinite subsets of $\mathbb{N}$. If there exist dense subsets 
		$\mathcal{B}_{1}, \mathcal{B}_{2}$ of $\mathcal{A}$ and a strictly increasing 
		sequence $\lbrace n_{k} \rbrace \subseteq \mathbb{N}$ such that for every fixed $p_{\alpha}$  it holds that either
		and 
		
		$$
		\lim_{k \to \infty} 
		\| \Phi_{n_{k}}^{-2} (b_{n_{k}}) \,
		\Phi_{n_{k}}^{-1} (b_{n_{k}}) \, y \|
		= 
		\lim_{k \to \infty} \| b^{-1}_{n_{k}} x \| 
		\, \| \Phi_{n_{k}}^{-1} (b_{n_{k}}) y \|
		=$$
		$$=
		\lim_{k \to \infty} \| \Phi_{n_{k}}^{-1}  (b_{n_{k}}) x \| 
		\, \| \Phi_{n_{k}}^{-1} (b_{n_{k}}) y \| = 0
		$$
		or
		$$\lim_{k \to \infty} \| \Phi_{n_{k}}  (b_{n_{k}}^{-1})b_{n_{k}}^{-1} y \| =
		\lim_{k \to \infty} \|  b_{n_{k}}^{-1} x ) \| \cdot \|   b_{n_{k}}^{-1} y \| 
		= \lim_{k \to \infty}  \|    \Phi_{n_{k}}^{-1}  (b_{n_{k}}) \| \|  b_{n_{k}}^{-1} y\|=0$$
		for all $x \in \mathcal{B}_{1}$ and $y \in \mathcal{B}_{2}$, 
		then the family  $\{ C_{n} \}_{n \in \mathbb{N}}$ is $\mathcal{F}$–semi-transitive.
	\end{corollary}
	\begin{remark}\label{Banach-bimodule}
		The above results and the proofs can be extended and generalized to the case when $ \mathcal{A} $ is a Banach bimodule over a unital normed algebra $ \mathcal{A}_1 .$ The condition \textit{(E)}  remains the same in this case, however, the condition $(\mathcal{P})$ should be then reformulated as follows:\\
		Let  $\{p_\alpha\}_\alpha$ be a set in $\mathcal{A} \cap \mathcal{A}_1 $ satisfying that given any open subset $O$ of $\mathcal{A}$ there exists some $ x \in \mathcal{A} \cup \mathcal{A}_1$ and some  $ p_{\alpha_0} \in \{p_\alpha\}_\alpha$ such that  $ p_{\alpha_0}x \in  \mathcal{A}_1$ and $ p_{\alpha_0}^{3}x \in  O.$
		
		Also, in this case, we let
		$
		\{\Phi_{t,1}\}_{t\in \mathcal{S}},\ \{\Phi_{t,2}\}_{t\in \mathcal{S}},\ \ldots,\ \{\Phi_{t,N}\}_{t\in \mathcal{S}}
		$
		be families of isometric algebra isomorphisms of $\mathcal{A}_1$ and $$
		\{\Psi_{t,1}\}_{t\in \mathcal{S}},\ \{\Psi_{t,2}\}_{t\in \mathcal{S}},\ \ldots,\ \{\Psi_{t,N}\}_{t\in \mathcal{S}}
		$$ be families of linear isometric isomorphisms of $ \mathcal{A} $ such that for all $ x, y \in \mathcal{A}_1 $ and $ a \in \mathcal{A} $ it holds that $ \Psi_{t,j}(xay) = \Phi_{t,j}(x)\Psi_{t,j}(a) \Phi_{t,j}(y) $   for all $t\in \mathcal{S} $ and $j \in \lbrace 1,\dots , N \rbrace.$ We express aperiodicity still in terms of the families\\
		$
		\{\Phi_{t,1}\}_{t\in \mathcal{S}},\ \{\Phi_{t,2}\}_{t\in \mathcal{S}},\ \ldots,\ \{\Phi_{t,N}\}_{t\in \mathcal{S}}
		,$ moreover, the conditions \textit{(R)} and \textit{(C)} remain the same in this case. Further, we define the operator
		$$
		T_{t,\ell} : \mathcal{A} \to \mathcal{A}
		\quad \text{by} \quad
		T_{t,\ell}(a) = b_{t,\ell} \, \Psi_{t,\ell}(a)
		$$  for all $ a \in \mathcal{A}.$ Under this set up and these assumptions, it is not hard to extend and generalize the proof of Theorem \ref{algebra} in this setting.
	\end{remark}
	
	\section{Applications}
	From now on and in the rest of this section we will assume that $\mathcal{F}$ is the family of all infinite subsets of $\mathbb{N}.$ Further, in this section, for a separable Hilbert space $H$ we let $B(H)$ be  
	the space of all bounded linear  operators on $H$, and  
	$\mathcal{C}$ be the $C^{*}$-algebra of compact operators on $H$. For an orthonormal basis $\{ e_{j} \}_{j \in \mathbb{Z}}$ for $H$,  
	we let for each $m \in \mathbb{N}$, $P_{m}$ be the orthogonal projection  
	onto $\text{Span} \{ e_{-m}, \dots, e_{m} \}$. 
	
	Moreover, we will denote by $\Omega$ be a locally compact, non-compact  Hausdorff, space. As usual, $C_{b}(\Omega)$  will denote the space of all bounded continuous functions on $\Omega$ and $C_{0}(\Omega)$ will denote the space of all continuous functions on $\Omega$ vanishing at infinity. Both these spaces will be equipped with supremum norm. In addition, we let $ C_{c}(\Omega) $ be the space of all continuous compactly supported functions on $\Omega .$ For each pair of  compact subsets $K_1$ and $K_2$ of $\Omega$ with $ K_1 \subseteq K_2 $ , we let $u_{(K_1, K_2)} \in C_{0}(\Omega)$  be 
	such that $u_{{(K_1, K_2)}|_{K_1}} = 1$, $\operatorname{supp} u_{(K_1, K_2)} = K_2$ and  $0 \leq u_{(K_1, K_2)} \leq 1$.
	
	\subsection{Applications to generalized weighted shifts on Hilbert C*-modules}
	In this subsection, we will denote by $\ell_{2}(\mathcal{C})$  
	the standard (right) Hilbert module over $\mathcal{C}$, see \cite[Example 1.3.5]{MT}. Notice that $\ell_{2}(\mathcal{C})$ is a Banach algebra. Indeed, we can define multiplication on  $\ell_{2}(\mathcal{C})$ as pointwise multiplication, i.e., if $\{x_j\}_{j\in\mathbb{Z}}, \{y_j\}_{j\in\mathbb{Z}} \in \ell_{2}(\mathcal{C})$,  
	then $$\{x_j \}_{j\in\mathbb{Z}} \cdot \{ y_j\}_{j\in\mathbb{Z}}=\{x_j y_j\}_{j\in\mathbb{Z}} .$$
	To see that $\{x_j y_j\}_{j\in\mathbb{Z}}$ belongs to 
	$\ell_{2}(\mathcal{C})$, 
	it suffices to observe that
	$$
	\sum_{j} y_j^{*} x_j^{*} x_j y_j 
	\le \sum_{j} y_j^{*} \|x_j\| ^{2} y_j 
	 \leq 
	\parallel  \{ x_j \}_{j \in \mathbb{Z}} \parallel^{2} \sum y_j^{*} y_j .
	$$
	
	For each $m, J \in \mathbb{N},$ we let 
	$\tilde{p}_{J,m} \in \ell_{2}(\mathcal{C})$ be given by
	$$
	(\tilde{p}_{J,m})_{i} =
	\begin{cases}
		P_{m}, & \text{if } -J \le i \le J, \\
		0, & \text{else}.
	\end{cases}
	$$
	
	As shown in the proof of \cite[Proposition 3.1]{FIL}, the set 
	$\{ \tilde{p}_{J,m} \}_{(J,m) \in \mathbb{N}^{2}}$
	forms a left approximate unit for 
	$ \ell_{2}(\mathcal{C})$. We put then $\mathcal{A} = \ell_{2}(\mathcal{C}),$ 
		and $$\mathcal{A}_1 := \{  \{ y_j \}_{j \in \mathbb{Z}} \mid y_j \in B(H) \ \forall j, \, \exists M_y > 0 
		\text{ such that } \| y_j \| \le M_y \ \forall j \} .$$
		The multiplication on $\mathcal{A}_1$ is defined similarly component-wise. We define the norm $ \| \cdot \|_1 $ on $\mathcal{A}_1$ as $ \| \{ y_j \}_{j \in \mathbb{Z}} \|_1 = sup_{j \in \mathbb{Z} } \| y_j \| $ for all $ \{ y_j \}_{j \in \mathbb{Z}} \in \mathcal{A}_1 .$ It is straightforward to check that the condition \textit{(E)} is satisfied in this case.
		
		 Let now $\mathcal{U}$ be a unitary operator on $H$. For each $ n \in \mathbb{N} $ 
		we define the map $\Phi_{n} : \mathcal{A}_1 \to \mathcal{A}_1$ by
		$$
		\Phi_{n} \left( \{ y_j \}_{j \in \mathbb{Z}} \right) 
		= \{ \mathcal{U}^{-n} y_{j - n}  \mathcal{U}^{n} \}_{j \in \mathbb{Z}},
		$$ for all $ \{ y_j \}_{j \in \mathbb{Z}} \in \mathcal{A}_1 ,$ 
		(that is the $j$-th coordinate of $\Phi_{n} (\{y_j\}_{j\in\mathbb{Z}}) = \mathcal{U}^{-n} y_{j - n}  \mathcal{U}^{n}$ for all $j \in \mathbb{Z}$).\\
		Since the set 
		$\{ \tilde{p}_{J,m} \}_{(J,m) \in \mathbb{N}^{2}}$
		forms a left approximate unit for 
		$\mathcal{A} =  \ell_{2}(\mathcal{C})$ and $\tilde{p}_{J,m}= \tilde{p}_{J,m}^{3} $ for all $ (J,m) \in \mathbb{N}^{2},$ it follows that the condition $ (\mathcal{P} )$ is satisfied in this case. Moreover, by some calculations it can be checked that $\{ \Phi_{n} \}_{n \in \mathbb{N}}$ 
		is aperiodic with respect to the left approximate unit $\{ \tilde{p}_{J,m} \}_{(J,m) \in \mathbb{N}^{2}} $ and that the condition \textit{ (R) } is satisfied in this case.  
		
		For each $n \in \mathbb{N},$ let $\{ W_{n,j} \}_{j \in \mathbb{Z}}$ be
		a family of operators in $B(H)$ which is uniformly bounded in norm such that each $W_{n,j}$ has a bounded inverse and such that $\{ W_{n,j}^{-1} \}_{j \in \mathbb{Z}}$ is also uniformly bounded in the norm. Then, for each $n \in \mathbb{N}$ 
		we let $b_{n} \in \mathcal{A}_1$ be given by
		$$
		(b_{n})_{i}
		= 
		W_{n,i}\mathcal{U}^{n},
		$$ for all $ i \in \mathbb{Z} .$ Obviously, $ b_{n} $   is invertible in $\mathcal{A}_1$ for each $n \in \mathbb{N}$ and 
		$$
		(b_{n}^{-1})_i=  
		\mathcal{U}^{-n} \, W^{\,-1}_{n, i},    
		$$ for all $ i \in \mathbb{Z} ,$
		hence, the condition $(C)$ is satisfied in this case.\\ 
		Finally, for each $n \in \mathbb{N}$  we let  $T_{n}$ be the operator on $\mathcal{A}$ given by
		$$
		T_{n}(a)
		= b_{n} \,
		\Phi_{n}(a),
		\qquad \text{for all } a \in \mathcal{A}.
		$$
		Recall now from the previous section that $ C_n = \frac{1}{2} (T_n + T_n^{-1}) $ for each $n \in \mathbb{N}.$
		 
		Noticing that for each $ a \in \mathcal{A}$ and $ j \in \mathbb{Z}$ it holds that $\|a_j\| \le \| a \|$, in view of Corollary \ref{cosinus-gust} and by arguments similar to those in the proofs of \cite[Theorem 3.2]{FIL} and \cite[Proposition 2.7]{BIMS}   it is not hard to deduce the following proposition in this special case.
		\begin{proposition}\label{cosinus-shift}
Under the above notation and assumptions, if there exist dense subsets 
			$H_{0}$ and $H_{1}$ of $H$ and a strictly increasing 
			sequence $\{ n_{k} \}_{k}$ of natural numbers 
			such that either
			
			$$
			\lim_{k \to \infty} 
			\| W^{-1}_{n_{k}, j} x \| \, \| W_{ n_{k}, j+n_{k} } y \|
			=
			\lim_{k \to \infty}
			\| W_{n_{k} , j + n_{k}} x \| \, 
			\| W_{n_{k} , j + n_{k}} y \|
			$$
			
			$$
			\lim_{k \to \infty}
			\| W_{n_{k} , j + 2 n_{k}} 	W_{n_{k} , j + n_{k}} y \|
			= 0,
			$$
			
			or
			
			$$
			\lim_{k \to \infty}
			\| W_{n_{k}, j + n_{k}} x \| \, 
			\| W^{-1}_{n_{k}, j} y \|
			=
			\lim_{k \to \infty}
			\| W^{-1}_{n_{k}, j} x \| \, 
			\| W^{-1}_{n_{k}, j} y \|=$$
			$$
			=
			\lim_{k \to \infty}
			\| W^{-1}_{n_{k}, j - n_{k}}
			W^{-1}_{n_{k}, j} y \|
			= 0,$$

			for all $x \in H_{0}$ and $y \in H_{1}$, then
			$\{ C_{n} \}_{n \in \mathbb{N}}$ is 
			$\mathcal{F}$–semi-transitive on $ \mathcal{A}$.
		\end{proposition}
		
		\begin{example}\label{primer-modul-kompakt}
			Let  $H = L^{2}(\mathbb{R})$. As in \cite[Example 3.3]{FIL}, choose some 
			$\{ w_j \}_{j \in \mathbb{Z}} \subseteq L^{\infty}(\mathbb{R})$ such that $w_j > 0$, 
			$w_{j}^{-1}  \in L^{\infty}(\mathbb{R})$ and $\parallel w_{j} \parallel ,\parallel w_{j}^{-1} \parallel  \leq M $ for all $j \in \mathbb{Z}$ and some $M>1$. Let $\{ c_j \}_{j \in \mathbb{Z}} \subseteq \mathbb{R}^{+}$ be such that $c_j \ge K$ for all 
			$j \in \mathbb{Z}$ and some $K>0$. For each $j \in \mathbb{Z},$ let $\alpha_j$ be the translation on 
			$\mathbb{R}$ given by $\alpha_j(t) = t - c_{j}$ for all $t \in \mathbb{R},$ and let then $ V_j \in B(H)$ be given by  $$V_{j}(f) = w_{j} \cdot(f \circ \alpha_j)$$ for all $ f \in H.$ For each $n \in \mathbb{N}, i \in \mathbb{Z}$, 
			let $W_{n, i}  \in B(H) $ be given by
			$$
			W_{n, i} 		:= 
			V_{i} V_{i-1} \cdots 
			V_{i-n+1}.
			$$ If there exists an $\varepsilon > 0$ such that  
			$w_{j} \chi_{[0,\infty)} \leq 1 - \varepsilon,w_{j} \chi_{(-\infty,0)} = 1$ for all $j \geq 0$ and $ w_j =1 $  
			for all $j < 0 $; then by the same arguments as in \cite[Example 3.3]{FIL} one can show that the conditions of Proposition \ref{cosinus-shift} are satisfied with  $H_{0} = H_{1} = C_{c}(\mathbb{R}) $ in this case. However, since $ w_j =1 $  
			for all $j < 0 ,$ it is not hard to check that the conditions of \cite[Proposition 3.5]{FIL} will not be satisfied. \\
			On the other hand, if there exists an $\varepsilon > 0$ such that
			$w_{j} \chi_{[0,\infty)} = 1, w_{j} \chi_{(-\infty,0)} \geq 1 + \varepsilon$ for all $j < 0 $ and $ w_j = 1$ for all $j \geq 0 ,$ then again the conditions of Proposition \ref{cosinus-shift} will be satisfied with  $H_{0} = H_{1} = C_{c}(\mathbb{R}) ,$ however, since  $ w_j = 1$ for all $j \geq 0 ,$ it can be checked that the conditions of \cite[Proposition 3.5]{FIL} will not be satisfied.
		\end{example}

		\text{ } 
		
		Next, for $N \in \mathbb{N},$ let $r_{1}, \dots, r_{N} \in \mathbb{N}$ with $r_{1} < r_{2} < \dots < r_{N}.$
		   
		Let $\alpha_{1}, \dots, \alpha_{N}$ be homeomorphisms of $\Omega$,  
		$\{ w_{j}^{(1)} \}_{j \in \mathbb{Z}}, \dots, \{ w_{j}^{(N)} \}_{j \in \mathbb{Z}} \subseteq C_{b}(\Omega)$.  
		We will assume that $w_{j}^{(\ell)} > 0$ for all $j \in \mathbb{Z}$ and $\ell \in \{1, \dots, N\}$ and that $\| w_{j}^{(\ell)} \|_{\infty} < M$  
		for all $j \in \mathbb{Z}$, $\ell \in \{1, \dots, N\}$  
		and some $M > 0.$ 
		Now we consider $\mathcal{A} = \ell_{2}(C_{0}(\Omega))$,  
		and $$\mathcal{A}_1 := \{  \{ y_j \}_{j \in \mathbb{Z}} \mid y_j \in C_{b}(\Omega) \ \forall j, \, \exists M_y > 0 
		\text{ such that } \| y_j \|_{\infty} \le M_y \ \forall j \} .$$ The multiplication and the norm on $ \mathcal{A}_1$ are defined in the same way as in the previous case. Once again, it is straightforward to check that the condition \textit{ (E) } is satisfied in this case.

		For $J \in \mathbb{N}$, we let $\tilde{p}_{(K_1, K_2, J)} \in \mathcal{A}$ be given by 
		$$(\tilde{p}_{(K_1, K_2, J)})_{j} =
		\begin{cases}
			u_{(K_1, K_2)}, & \text{if } -J \leq j \leq J, \\[4pt]
			0, & \text{else.}
		\end{cases}$$

		Then it is easily checked that  
		$\{ \tilde{p}_{(K_1, K_2, J)}^{3} \}_{K_1 \subseteq K_2 \subseteq \Omega, \, K_1, K_2 \text{ compact ,} J \in \mathbb{N}}$  
		satisfies the condition  $ (\mathcal{P}.)$\\
		  For each $\ell \in \{1, \dots, N\}$ and $k \in \mathbb{N}$, we let $\tilde{\Phi}_{k, \ell} : \mathcal{A} \to \mathcal{A}$ be given by  
		$$
		\tilde{\Phi}_{k,\ell}\big( \{ f_{j} \}_{j \in \mathbb{Z}} \big)
		= 
		\{ f_{j-k} \circ \alpha_{\ell}^{k} \, \}_{j \in \mathbb{Z}},
		$$
		for all $\{ f_{j} \}_{j \in \mathbb{Z}} \subset \mathcal{A}$.  Let now $r_{1}, \dots, r_{N} \in \mathbb{N}$ with $r_{1} < r_{2} < \dots < r_{N}$.  
		For each $n \in \mathbb{N}$ and $\ell \in \{1, \dots, N\}$, set  
		$$
		\Phi_{n,\ell} = \tilde{\Phi}_{r_{\ell} n, \ell}.
		$$
		Clearly, the system of isometric isomorphisms $\{ \Phi_{k,1} \}_{k \in \mathbb{N}}, \dots, \{ \Phi_{k,N} \}_{k \in \mathbb{N}}$ is  
		disjoint aperiodic with respect to $$\{ \tilde{p}_{(K_1, K_2, J)} \}_{K_1 \subseteq K_2 \subseteq \Omega, \, K_1, K_2 \text{ compact ,} J \in \mathbb{N}}.$$   
		Also, the condition \textit{(R) } is satisfied in this case.\\
		 For each $ n \in \mathbb{N}$ and $\ell \in \{1, \dots, N\}$, let $ b_{n,\ell} \in \mathcal{A}_1 $ be given by  
		$$
		(b_{n,\ell})_{j+r_{\ell} n} = 
		\prod_{k=1}^{r_{\ell} n}
		\big( w_{j+k}^{(\ell)} \circ \alpha_{\ell}^{n-k} \big)
		\quad \text{for all } j \in \mathbb{Z}.
		$$ Since for each $ n \in \mathbb{N}$ and $\ell \in \{1, \dots, N\},$ the function $w_{j}^{(\ell)} > 0,$ it follows that the condition $(C)$ is satisfied in this case.\\
		Now, for each $ K_2 \subseteq \Omega $ compact, we can find some $ v_{K_2  } \in C_{0}(\Omega)$ with $  0 \leq v_{K_2  } \leq 1 $ such that ${v_{K_2}}_{|_{K_2}} = 1 .$ Recalling that the approximate unit for $ \mathcal{A}$ in this case is indexed over all triples $(K_1, K_2, J) $ with $ K_1 \subseteq K_2 \subseteq \Omega, \, K_1, K_2 \text{ compact, } J \in \mathbb{N} ,$ for each $ n \in \mathbb{N}$ and $\ell \in \{1, \dots, N\}$, and every fixed triple $ (K_1, K_2, J) $ we let $ b_{n,\ell, (K_1, K_2, J)}^{-1} \in \mathcal{A}_1 $ be given by $$ (b_{n,\ell, (K_1, K_2, J)}^{-1})_{j+r_{\ell} n} =\frac{v_{K_2}}{\prod_{k=1}^{r_{\ell} n}
			\big( w_{j+k}^{(\ell)} \circ \alpha_{\ell}^{n-k} \big) } $$ for all $ j \in \mathbb{Z} $ with $ -J - r_{\ell} n \leq j \leq J - r_{\ell} n$ and $(b_{n,\ell, (K_1, K_2, J)}^{-1})_{j+r_{\ell} n} = 0 $ else.
		Then, for all $n \in \mathbb{N}$ and $\ell \in \{1, \dots, N\},$ we let $T_{n,\ell} : \mathcal{A} \to \mathcal{A}$ be given by  
		$$
		T_{n,\ell}(a) = b_{n,\ell} \, \Phi_{n,\ell}(a)
		\quad \text{for all } a \in \mathcal{A}.
		$$ 
		  
		The next corollary can be deduced by some calculations from Theorem \ref{algebra}.
		
		\begin{corollary}\label{modul}
			Under the above notation and assumptions,  
			the following statements are equivalent:  
			
			(1) The sequences of operators  
			$$
			\{ T_{n,1} \}_{n \in \mathbb{N}}, \ \dots , \ \{ T_{n, N} \}_{n \in \mathbb{N}}
			$$  
			are $d\mathcal{F}$-semi-transitive on $\mathcal{A}$.  
			
			(2) For every $J \in \mathbb{N}$ and each compact subset $K$ of $\Omega$,  
			there exists a strictly increasing sequence $\{ n_{k} \}_{k} \subseteq \mathbb{N}$ such that  
			$$
			\lim_{k \to \infty} 
			\Bigg(
			\sup_{t \in K}
			\frac{1}{\prod_{i = 1}^{r_{s} n_{k}}
				\big( w_{j + i - r_{s} n_{k}}^{(s)} \circ \alpha_{s}^{r_{s} n_{k}-i} \big)(t)}
			\Bigg) \cdot
			\Bigg(
			\sup_{t \in K}
			\prod_{i = 1}^{r_{\ell} n_{k}}
			\big( w_{j + i}^{(\ell)} \circ \alpha_{\ell}^{-i} \big)(t)
			\Bigg)
			= 0
			$$
			for all $j \in \mathbb{Z}$ with $-J \le j \le J$ and every $l,s \in \{1, \dots, N\}$.  
			Moreover, for each distinct $s, \ell \in \{1, \dots, N\}$ and all $j \in \mathbb{Z}$  
			with $-J \le j \le J$, it holds that  
			$$
			\lim_{k \to \infty}
			\Bigg(
			\sup_{t \in K}
			\frac{
				\prod_{i = 1}^{r_{\ell} n_{k}}
				\big( w_{j - r_{\ell} n_{k} - i}^{(\ell)} 
				\circ \alpha_{\ell}^{-i } \circ \alpha_{s}^{r_{s} n_{k}} \big)(t)
			}{
				\prod_{i = 1}^{r_{s} n_{k}}
				\big( w_{j -r_{s} n_{k} + i}^{(s)} 
				\circ \alpha_{s}^{n_{k}-i} \big)(t)
			}
			\Bigg)
			= 0.
			$$
		\end{corollary}
		
		\begin{example}\label{primer-shift}
			Let $N = 1, r_1 =1$ and  choose some 
			$\{ w_j \}_{j \in \mathbb{Z}} \subseteq L^{\infty}(\mathbb{R})$ such that $w_j > 0$, 
			$w_{j}^{-1}  \in L^{\infty}(\mathbb{R})$ and $\parallel w_{j} \parallel \leq M $ for all $j \in \mathbb{Z}$ and some $M>1$. Let $ \alpha(t)= t-1 $ for all $ t \in \mathbb{R}.$  If there exists an $\varepsilon > 0$ such that  
			$w_{j} \chi_{[0,\infty)} \leq 1 - \varepsilon,w_{j} \chi_{(-\infty,0)} = 1$ for all $j \geq 0$ and $ w_j =1 $  
			for all $j < 0 $; then the conditions of Corollary \ref{modul} are satisfied in this case. However, since $ w_j =1 $  
			for all $j < 0 $; it follows that the conditions of \cite[Theorem 2.2]{CAOT} and \cite[Corollary 2.4]{CAOT} will not be satisfied in this case. \\
			On the other hand, if there exists an $\varepsilon > 0$ such that
			$w_{j} \chi_{[0,\infty)} = 1, w_{j} \chi_{(-\infty,0)} \geq 1 + \varepsilon$ for all $j < 0 $ and $ w_j = 1$ for all $j \geq 0 ,$ then the conditions of Corollary \ref{modul} are again satisfied, however, since  $ w_j = 1$ for all $j \geq 0 ,$ the conditions \cite[Theorem 2.2]{CAOT} and \cite[Corollary 2.4]{CAOT} will not be satisfied in this case.
		\end{example}
		
 \subsection{Applications to generalized weighted composition operators on spaces of operator-valued continuous functions}
 
 \text{ }		
 
  In this subsection, we let  
		$C_{0}(\Omega, \mathcal{C})$ denote the non-unital normed algebra of all continuous $\mathcal{C}-$valued functions 
		on $ \Omega$  vanishing at infinity and equipped with the norm $$\| f \|_{\infty} := \sup_{x \in \Omega} \| f(x) \| ,$$ 
		and $C_{b}(\Omega, B(H))$ be the unital normed algebra consisting of all uniformly bounded $B(H)-$valued continuous 
		functions on $\Omega$ equipped with the same norm $\| \cdot \|_{\infty}$.\\
		Obviously, $C_{0}(\Omega, \mathcal{C})$ is (left) ideal in $C_{b}(\Omega, B(H))$
		and the condition $(E)$ is satisfied. We will
		therefore denote $C_{0}(\Omega,  \mathcal{C})$ and 
		$C_{b}(\Omega, B(H))$ by $ \mathcal{A}$
		and $A_{1},$ respectively. For every  triple $(K_1, K_2, m) $ with $ K_1 \subseteq K_2 \subseteq \Omega, \, K_1, K_2 \text{ compact, } m \in \mathbb{N} ,$ let $\tilde{P}_{(K_{1},K_{2},m)} \in C_{0}(\Omega, \mathcal{C})$ be given by
		\[
		\tilde{P}_{(K_{1},K_{2},m)}(x) = u_{(K_{1},K_{2})}(x) P_{m} \quad \text{for all } x \in \Omega .
		\]
		We have the following lemma.
		
		\begin{lemma}\label{cover}
			The set 
			$
			\lbrace \tilde{P}_{(K_{1},K_{2},m)} \rbrace_{K_{1} \subseteq K_{2} \subseteq \Omega, K_{1}, K_{2} \text{ compact, } 
				m \in \mathbb{N} }
			$
			satisfies the condition $(\mathcal{P})$.
		\end{lemma}

		\begin{proof}
			Let $O$ be a non-empty open subset of $ \mathcal{A}$ and 
			$f \in O$. Then there exists some $\varepsilon > 0$ such that $\varepsilon$–
			neighbourhood of $f$ is contained in $O$ as well. Since 
			$f \in \mathcal{A},$ there exists some compact subset $K_{1}$ of $ \Omega$ such 
			that $\| f(x) \| < \frac{\varepsilon}{2}$ for all $x \in \Omega \setminus K_{1}$. 
			Now, by continuity of $f,$ for each $x \in K_{1}$ we can find some open set $O_{x}$ of $\Omega$ containing $x$ such that $\| f(u) - f(x) \| < \frac{\varepsilon}{2}$ whenever 
			$y \in O_{x}$. Then for all $m \in \mathbb{N}$ and each $y \in O_{x}$ we obtain
			
			\[
			\| (I - P_{m}) f(y) \| \le \| (I - P_{m}) ( f(y) - f(x) ) \| 
			+ \| (I - P_{m}) f(x) \|
			\]
			
			\[
			\le \| f(y) - f(x) \| + \| (I - P_{m}) f(x) \| 
			< \frac{\varepsilon}{2} + \| (I - P_{m}) f(x) \|
			\]

			By \cite[Proposition 2.2.1]{MT}, since $f(x) \in \mathcal{C},$ we can find sufficiently large $m_{x} \in \mathbb{N}$ such that 
			$ \| (I - P_{m}) f(x) \| < \frac{\varepsilon}{2} $ for all $ m \geq m_x.$ 
			Hence, $\|  (I - P_{m}) f(y) \| < \varepsilon $ each $y \in O_{x}$ and all $ m \geq m_x.$ 
			
			The collection 
			$
			\{O_{x} \}_{x \in K_{1}} $ is an open cover of $K_{1}$, hence there exist a 
			finite subcover $ \{ O_{x_{1}}, \ldots, O_{x_{M}} \}.$
			Let 
			$\tilde{m} := \max \{ m_{x_{1}}, \ldots, m_{x_{M}} \}. $
			Then, by the above we get

			$ \| \tilde{P}^{3}_{(K_{1},K_{2},m)}(x) f(x) - f(x) \| 
			= \| (I - P_{m}) f(x) \| < \varepsilon$
			for all $ x \in K_{1}$ and $ m \geq \tilde{m}.$ 
			If $ x \in \Omega \setminus K_{1},$  then for all $ m \in \mathbb{N} $ we get
			\[
			\| \tilde{P}^{3}_{(K_{1},K_{2},m)} (x)f(x) \| 
			+ \| f(x) \| 
			= \| u^{3}_{(K_{1},K_{2})}(x) P_{m} f(x) \|
			+ \| f(x) \| 
			\le 2 \| f(x) \| < 2 \frac{\varepsilon}{2} = \varepsilon.
			\]
			
			Thus,
			\[
			\| \tilde{P}^{3}_{(K_{1},K_{2},\tilde{m})} f - f \|_{\infty} < \varepsilon,
			\text{ so } \tilde{P}^{3}_{(K_{1},K_{2},m)} f \in O. 
			\]
		\end{proof}

		Now we consider families $\{ \alpha_{n,1} \}_{n \in \mathbb{N}} \dots,\{ \alpha_{n,N} \}_{n \in \mathbb{N}}$ of homeomorphisms of $ \Omega .$  and families of unitary operators $\{ U_{n,1} \}_{n \in \mathbb{N}} \dots,\{ U_{n,N} \}_{n \in \mathbb{N}}$ on $ H.$
		For every $n \in \mathbb{N}$ and $l \in \{1, \ldots, N \}$ we let then
		
		\[
		\Phi_{n,l} :\mathcal{A}_{1} \rightarrow \mathcal{A}_{1}
		\text{ be given by } 
		\Phi_{t,l}(f) := U_{n,l}^{*} (f \circ \alpha_{n,l})U_{n,l} 
		\]
		for all $ f \in \mathcal{A} .$
		Obviously the system $\{ \Phi_{n,1} \}_{n \in \mathbb{N}} , \dots , \{\Phi_{n,N} \}_{n \in \mathbb{N}} $ 
		is a system of isometric isomorphisms on $ \mathcal{A}$ which 
		is disjoint aperiodic with respect to 
		$\{ \tilde{P}_{(K_{1},K_{2},m)}  \}$ and,  
		clearly, the condition ($R$) is satisfied. 
		Let $\{ w_{n,1} \}_{n \in \mathbb{N}} , \dots , \{ w_{n,N} \}_{n \in \mathbb{N}} $ be families in $\mathcal{A}_{1}$ such that 
		$w_{n,l}(x)$ is invertible for all $x \in \Omega$, 
		$n \in \mathbb{N}$, $l \in \{1, \ldots, N\}$, and such that 
		$\{ w^{-1}_{n,1} \}_{n \in \mathbb{N} } , \dots ,\{ w^{-1}_{n,N} \}_{n \in \mathbb{N} } $ also belong to $\mathcal{A}_{1}$ i.e. 
		
		\[
		\sup_{x \in \Omega} \| w^{-1}_{n,l}(x) \| < \infty 
		\quad \text{for all } n \in \mathbb{N} \text{ and } l \in \{1, \ldots, N\}.
		\]

		For each $ n \in \mathbb{N}$, and $l \in \{1, \ldots, N\},$ let then $b_{n,l} := w_{n,l}U_{n,l}$.
		Since each $b_{n,l}$ is invertible in $\mathcal{A}_{1}$, it follows that
		the condition $(C)$ is satisfied in this case.
		\begin{corollary}\label{operator-algebra}
			Under the above notation and assumptions, if there exist dense 
			subsets $H_{0}, H_{1}, \ldots, H_{N }$ of $H$ such that for each 
			$s, l \in \{1, \ldots, N\} $ and every compact subset $K$ of $\Omega$ we have that

			\[
			\lim_{n \to \infty} 
			\left[ 
			\left( \sup_{x \in K} \| w_{n,s}^{-1}(x) h_{s} \| \right) \cdot
			\left( \sup_{y \in K} \| (w_{n,\ell} \circ \alpha_{n,\ell}^{-1})(y) h_{0} \| \right)
			\right]
			= 0
			\]
			
			for every $h_{s} \in H_{s}$ and $h_{0} \in H_{0},$ and, in addition,
			for each distinct $s, l \in \{1, \ldots, N\}$ and every $h_{s} \in H_{s}$ we have
			
			\[
			\lim_{n \to \infty} \left(  \sup_{x \in K}
			\| (w_{n,l} \circ \alpha_{n,l}^{-1} \circ \alpha_{n,s})   (x) \, w_{n,s}^{-1}(x) \, h_{s} \|  \right) = 0,
			\]
			
			then the families $\{ T_{n,1} \}_{n \in \mathbb{N}}, \ldots, 
			\{ T_{n,N} \}_{n \in \mathbb{N}}$ are $d\mathcal{F}$-semi-transitive.
		\end{corollary}

		\begin{proof}
			
			For each $s \in \{1, \ldots, N\}$, the set 
			consisting of elements of the form $$u_{(K_{1},K_{2})} G_{s,m} ,$$
			where $G_{s,m} \in \mathcal{C}$ is such that
			\[
			G_{s,m} e_{i} \in H_{s} \text{ for all } j \in \{ -m, \ldots, m \}  \text{ and } G_{s,m}e_{j} = 0 \text{ for } j \in \mathbb{Z} \setminus \{-m, \ldots, m \}
			\]
			is obviously dense in $ \mathcal{A}$ by Lemma \ref{cover} and the density of $H_{s}$ in $H$. (Here of course $m $ runs through the whole $\mathbb{N} $ and $ (K_{1},K_{2}) $ runs through all pairs of compact subsets $ K_{1},K_{2}$ of $ \Omega $ with $ K_{1} \subseteq K_{2}.$) 
			Let us denote this set by $\mathcal{A}_{s}$ for each $s \in \{1, \ldots, N\}$. Similarly, we let $ \mathcal{A}_{0}$ be the set consisting of all elements of 
			the form $u_{(K_{1},K_{2})} D_{m}$ where $D_{m} \in \mathcal{C}$ is such that
			
			\[
			D_{m} e_{j}\in H_{0} \text{ for all } j \in \{-m, \dots m\}
			\quad \text{and} \quad
			D_{m} e_{j} = 0 \text{ for } j \in \mathbb{Z} \setminus \{-m, \dots m\}.
			\]
			
			By some calculations similar to those in the proof of \cite[Theorem 3.2]{FIL} and \cite[Proposition 2.7]{BIMS}, it is then not hard to deduce that 
			the conditions of Corollary \ref{algebra-primena} are satisfied.
		\end{proof}

		\begin{example}\label{primer-operator-algebra}
			Let $\Omega = \mathbb{R} \text{ and } H = L^{2}(\mathbb{R}). $
			Put $N = 1$ and let $ \eta \in C_{b}(\mathbb{R}^{2})$ be such that 
			$\eta(x,t) = 1 \text{ whenever } x \le 0 \text{ or } t \le 0 
			\text{ and that } \eta(x,t) = \tfrac{1}{2} \text{ whenever } x, t \ge 1.$
			Assume also that $\tfrac{1}{2} \le \eta(x,t) \le 1 \text{ for all } (x,t) \in \mathbb{R}^{2}.$
			For each $n \in \mathbb{N},$ let $\alpha_{n} : \mathbb{R} \to \mathbb{R} $
			be given by $\alpha_{n}(t) = t - n \text{ for all } t \in \mathbb{R}.$
			For each $x \in \mathbb{R} ,$ let $ W_{x} \in B(H)$ be given 
			$W_{x}(f) = \eta(x,\cdot) \, ( f \circ \alpha_{1}  ) \text{ for all } f \in H$
			(that is $W_{x}(f)(t) = \eta(x,t) \cdot f(t-1)$ for all $t \in \mathbb{R}$ and $f \in H$).
			For each $n \in \mathbb{N},$ let then $w_{n} \in A_{1}$ be given by 
			$w_{n}(x) = W_{x} W_{x-1} \dots W_{x-n}$ and let $H_{0} = H_{1}=C_{c}(\mathbb{R}).$ By some calculations it can be checked that the conditions of Corollary \ref{operator-algebra} are satisfied in this case.
		\end{example}

			\section{ Disjoint $\mathcal{F}-$semi-transitivity in weighted solid spaces}
			In this section, $X$ will denote a topological space and $ \alpha $ will denote a homeomorphism on $ X ,$ 
			We recall also the following definitions.
			%=============================================================
			\begin{definition}
				A Banach function space $\mathcal{Y}$ on $X$ is called \emph{solid} if for each $f\in \mathcal{Y}$ and $g\in\mathcal M_0(X)$, satisfying $|g|\leq |f|$, we have $g \in \mathcal{Y}$
				and $\|g\|_{\mathcal{Y}}\leq \|f\|_{\mathcal{Y}}$.
			\end{definition}
			For the next results we shall also assume that the following conditions from \cite{saw} on the Banach function space $\mathcal{Y}$ hold.
			
			\begin{definition}\label{condition}
				Let $X$ be a topological space, $\mathcal{Y}$ be a Banach function space on $X$, and $\alpha$ be an aperiodic homeomorphism of $X .$ We say that $\mathcal{Y}$ satisfies condition $\Omega_\alpha$ if the following conditions hold:
				~\	\begin{enumerate}
					\item $\mathcal{Y}$ is solid and $\alpha$-invariant;
					\item for each compact set	$K\subseteq X$ we have $\chi_{K}\in \mathcal{Y}$;
					\item $\mathcal{Y}_{bc}$ is dense in $\mathcal{Y}$, where $\mathcal{Y}_{bc}$ is the set of all bounded compactly supported functions in $\mathcal{Y}$.
				\end{enumerate}
			\end{definition}
			Moreover, we let $\mathcal{Y}_{b}$ be denote the unital algebra of all bounded measurable functions on $X$ equipped with the supremum norm.
			
			From now on, we shall assume that $\mathcal{Y}$ is a Banach function space satisfying the conditions of all the definitions above. For a measurable positive function $ w $ on $X ,$ we let $ w^{-1} := \frac{1}{w}  .$ If $ w  $ is a positive measurable function on $  X $ such that $ w $ is bounded, then $ T_{\alpha,w} $ will denote the weighted composition operator on $ \mathcal{Y}  $ defined by $ T_{\alpha,w} (f) = w \cdot (f \circ \alpha)  $ for all $ f \in  \mathcal{Y}.  $  
			Let now $ \{ \alpha_{t,1} \}_{t \in S} \dots,\{ \alpha_{t,N} \}_{t \in S} $ be families of homeomorphisms of $X$ and $$ \{ w_{t,1} \}_{t \in S} \dots,\{ w_{t,N} \}_{t \in S} $$ be families of bounded positive measurable functions on $X.$  We will assume that the system $\{ \alpha_{t,1} \}_{t \in S} \dots,\{ \alpha_{t,N} \}_{t \in S}$ is disjoint aperiodic, that is for every compact subset $K$ of $X$ and every $H \in \mathcal{F}$  
			there exists some $F \in \mathcal{F}$ with $F \subseteq H$ such that  
			
			\[
			K \cap \alpha_{t,\ell}(K) = \varnothing \quad \text{and} \quad
			\alpha_{t,\ell}^{-1} (\alpha_{t,s}(K)) \cap K = \varnothing 
			\]
			
			for all $t \in F$ and each distinct $s, \ell \in \{ 1,\dots, N \},$  and we will assume  that $\mathcal{Y}$ is $\alpha_{t,\ell}$-invariant for each $ t \in \mathcal{S}, \ell \in \{ 1,\dots, N \} .$ Also, we will suppose that $ w_{t,\ell}^{-1} $ is bounded on compact subsets of $X$ for each $ t \in \mathcal{S}, \ell \in \{ 1,\dots, N \} .$  
			
			For each $ t \in \mathcal{S}, \ell \in \{ 1,\dots, N \} $ we let then $T_{t,\ell} $ be the corresponding weighted composition operator on $\mathcal{Y}$ given by $T_{t,\ell}(f) =  w_{t,\ell} \cdot ( f \circ \alpha_{t,\ell} ) $ for all $ f \in \mathcal{Y} .$ 
			
			It is not hard to check that the conditions of Remark \ref{Banach-bimodule} are satisfied with $ \mathcal{A} = \mathcal{Y}, \mathcal{A}_1 = \mathcal{Y}_{b},  \{p_\alpha\}_\alpha= \{\chi_K\}_{K \subseteq X, K compact} $ and $$ \Psi_{t,\ell}(f) = f \circ \alpha_{t,\ell} , \text{ } \Phi_{t,\ell}(g) = g \circ \alpha_{t,\ell}$$ for all $ f \in \mathcal{Y}, g \in \mathcal{Y}_{b} $ and  each $ t \in \mathcal{S}, \ell \in \{ 1,\dots, N \} .$  Further, for each $ t \in \mathcal{S}, \ell \in \{ 1,\dots, N \} $ we let $ b_{t,\ell} := w_{t,\ell} .$ Then, for every compact subset $K$ of $X$ we have $$b_{t,\ell, K}^{-1} := \chi_{K} w_{t,\ell}^{-1} .$$ It follows then that $T_{t,\ell}(f) = b_{t,\ell} \, \cdot \Psi_{t,\ell}(f)
			$  for all $ f \in \mathcal{A} $ and each $ t \in \mathcal{S}, \ell \in \{ 1,\dots, N \} .$ Therefore, thanks to Theorem \ref{algebra}, Remark \ref{Banach-bimodule} and some calculations, we can deduce the following theorem.
			
			\begin{theorem}\label{tabatabaie}
				Under the above notation and assumptions, the following statements are equivalent.
				\begin{enumerate}
					\item The families $\{ T_{t,1} \}_{t \in \mathcal{S}}, \cdots,\{ T_{t,N} \}_{t \in \mathcal{S}}$ are $d \mathcal{F}$-semi-transitive on $\mathcal{Y}.$
					\item For each $\varepsilon \in (0,1)$ and every compact subset $K $ of $ X$  
					there exist some $F \in \mathcal{F}$ and a family 
					$\{ E_{t} \}_{t \in F}$  
					of Borel subsets of $K$ such that  
					
					\[
					\| \chi_{\, K \setminus E_{t}} \|_{\mathcal{Y}} < (4 + 2N)N\,\varepsilon ,
					\]
					
					\[
					\Bigl( \sup_{x \in E_t}
					\Bigl(  \,
					w_{t,s}^{-1}(x)  \Bigr) \Bigr) \cdot
					\;
					\Bigl( \sup_{x \in E_t} \Bigl(
					(w_{t,\ell} \circ \alpha_{t,\ell}^{-1})(x) \Bigr)
					\Bigr)
					<
					\frac{m_{K}^{2}\, \varepsilon^{2}}{(1-\varepsilon)^{2}}
					\]
					
					for all \( t \in F \) and \( \ell, s \in \{1,\dots,N\} \), and moreover  
					for all \( t \in F \) and each distinct \( s,\ell \in \{1,\dots,N\} \) it holds that
					
					\[
					\sup_{x \in E_t}
					\Biggl[
					\frac{
						(w_{t,\ell} \circ \alpha_{t,\ell}^{-1} \circ \alpha_{t,s})(x)
					}{
						w_{t,s}(x)
					}
					\Biggr]
					<
					\frac{m_{K}\, \varepsilon}{1 - \varepsilon}.
					\]
					
				\end{enumerate}
			\end{theorem}

			A weight on $X$ with respect to a homeomorphism $\alpha$ of $X$ is a continuous function $\eta : X \to (0, \infty)$ which satisfies
			$$
			\eta(\alpha(x)) \leq K_{\alpha} \eta(x) \quad \text{and} \quad \eta({\alpha}^{-1}(x)) \leq K_{\alpha} \eta(x) \quad \text{for all } x \in X,
			$$
			and some constant $K_{\alpha}
			> 0$.\\
			Below is an example of a such function.

			\begin{example}\label{primer-tezina} Let $X = \mathbb{R}$ and $p > 1$. Put
				$$
				\eta(x) = 
				\begin{cases}
					1\;\;\;\;\;\;\; \text{for} & x \in [-1, \infty), \\
					\frac{1}{|x|^{p}}\;\;\;  \text{for} & x \in (-\infty, -1).
				\end{cases}
				$$
				Then $\eta$ is a weight on $\mathbb{R}$ with respect to the translation $\alpha$ on $\mathbb{R}$ given by
				$$
				\alpha(x) = x + 1, \quad \text{for all } x \in \mathbb{R}.
				$$ Similar conclusion holds if we instead consider  $\alpha(x) = x - 1$, for all $x \in \mathbb{R}$.
			\end{example}
			
			 We define the weighted solid Banach function space $\mathcal{Y}_{\eta}$ on $X$ by
			$$
			\mathcal{Y}_{\eta} = \left\{ f \in \mathcal{M}_{0}(X) \ \middle| \ f \eta \in \mathcal{Y} \right\}
			,$$
			and we equip $\mathcal{Y}_{\eta}$ with the norm $\|\cdot\|_{\mathcal{Y}_{\eta}}$ given by
			$$
			\|f\|_{\mathcal{Y}_{\eta}} = \|f \eta\|_{\mathcal{Y}} \quad \text{for all } f \in \mathcal{Y}_{\eta}.
			$$

			\begin{lemma}\cite[Lemma 2.5]{arxiv2}
				The operator $T_{\alpha,w}$ is a bounded, linear self-mapping on $\mathcal{Y}_{\eta}$.
			\end{lemma} 
			
			In what follows, we will assume that $ \eta$ is a weight with respect to   $\alpha_{t,\ell}$ for each $ t \in \mathcal{S}, \ell \in \{ 1,\dots, N \} .$ We notice that, while $\|\cdot\|_{\mathcal{Y}}$ is $\alpha_{t,\ell}-$invariant, the norm $\|\cdot\|_{\mathcal{Y}_{\eta}}$ does not need to be $\alpha_{t,\ell}-$invariant.

			\begin{theorem}\label{vekt}
				Under the above notation and assumptions, the following statements are equivalent.
				\begin{enumerate}
					\item The families $\{ T_{t,1} \}_{t \in \mathcal{S}}, \cdots,\{ T_{t,N} \}_{t \in \mathcal{S}}$ are $d \mathcal{F}$-semi-transitive on $\mathcal{Y}_{\eta}.$
					\item For each $\varepsilon \in (0,1)$ and every compact subset $K $ of $ X$  
					there exist some $F \in \mathcal{F}$ and a family 
					$\{ E_{t} \}_{t \in F}$  
					of Borel subsets of $K$ such that  
					
					\[
					\| \chi_{\, K \setminus E_{t}} \|_{\mathcal{Y}} < (4 + 2N)N\,\varepsilon ,
					\]
					
					\[
					\Bigl( \sup_{x \in E_t}
					\Bigl( (\eta \circ \alpha_{t,s})(x) \,
					w_{t,s}^{-1}(x)  \Bigr) \Bigr) \cdot
					\;
					\Bigl( \sup_{x \in E_t} \Bigl(
					(\eta \circ \alpha_{t,\ell}^{-1})(x)
					(w_{t,\ell} \circ \alpha_{t,\ell}^{-1})(x) \Bigr)
					\Bigr)
					<
					\frac{m_{K}^{2}\, \varepsilon^{2}}{(1-\varepsilon)^{2}}
					\]
					
					for all \( t \in F \) and \( \ell, s \in \{1,\dots,N\} \), and moreover  
					for all \( t \in F \) and each distinct \( s,\ell \in \{1,\dots,N\} \) it holds that
					
					\[
					\sup_{x \in E_t}
					\Biggl[
					(\eta \circ \alpha_{t,\ell}^{-1} \circ \alpha_{t,s})(x)
					\frac{
						(w_{t,\ell} \circ \alpha_{t,\ell}^{-1} \circ \alpha_{t,s})(x)
					}{
						w_{t,s}(x)
					}
					\Biggr]
					<
					\frac{m_{K}\, \varepsilon}{1 - \varepsilon}.
					\]
					
				\end{enumerate}
			\end{theorem} 
			
			\begin{proof} 
				We prove first	(1) $\Rightarrow$ (2).  
				Let $K \subseteq X$, $K$ compact and $\varepsilon \in (0,1)$ be given. Put
				$$
				m_K = \inf_{x \in K} \eta(x).
				$$ and let  $ B\bigl( \chi_{K} , m_K \varepsilon^{2}\bigr)$ denote the open ball in $\mathcal{Y}_\eta $  with centre in $ \chi_{K}$ and radius $ m_K \varepsilon^{2}.$ 
				Choose some $H \in \mathcal{F}$ such that for all $t \in H$
				there exists some $\lambda_t \in \mathbb{R}^{+}$ satisfying  
				
				\[
				B\bigl( \chi_{K} , m_K \varepsilon^{2}\bigr) \;\cap\;
				\frac{1}{\lambda_{t}}T_{t,1}^{-1} \bigl( B(\chi_{K}, m_k \varepsilon^{2}) \bigr) \cap \cdots \cap \frac{1}{\lambda_{t}}T_{t,N}^{-1} \bigl( B(\chi_{K}, m_k \varepsilon^{2}) \bigr) \neq \varnothing .
				\]
				
				Since the system $\{ \alpha_{t,1} \}_{t \in S} \dots,\{ \alpha_{t,N} \}_{t \in S}$ is disjoint aperiodic,  
				there exists some $F \in \mathcal{F}$ with $F \subseteq H$ such that  
				
				\[
				K \cap \alpha_{t,\ell}(K) = \varnothing \quad \text{and} \quad
				\alpha_{t,\ell}^{-1} (\alpha_{t,s}(K)) \cap K = \varnothing 
				\]
				
				for all $t \in F$ and each distinct $s, \ell \in \{ 1,\dots, N \}.$  
				Since 
				
				\[
				\frac{1}{\lambda_{t}}   T_{t,\ell}^{-1} \bigl( B( \chi_{K}, m_K \varepsilon^{2}) \bigr)
				\;\cap\; B\bigl( \chi_{K}, m_K \varepsilon^{2} \bigr) \neq \varnothing ,
				\]
				for each $t \in F$ and each $\ell \in \{1,\dots, N\}$  
				we can find some $\widetilde{f}_{t,\ell} \in \mathcal{Y}_{\eta}$ such that  
				
				\[
				\|\widetilde{f}_{t,\ell} - \chi_{K}\|_{\mathcal{Y}_{\eta} } < m_k \varepsilon^{2}
				\]
				and
				
				\[
				\| \lambda_{t} T_{t,\ell} (\widetilde{f}_{t,\ell} - \chi_{K}) \|_{\mathcal{Y}_{\eta} } < m_k \varepsilon^{2}.
				\]
				
				Put  
				\[
				f_{t,\ell} :=\lambda_{t} T_{t,l} (\tilde{f}_{t,\ell})
				%( w_{t,\ell}  \circ    \alpha_{t,\ell}^{-1}  )^{-1}  ( f_{t,\ell}  \circ    \alpha_{t,\ell}^{-1}  )^{-1}  \frac{1}{\lambda_{t}} .
				\]
				then  
				
				\[
				\tilde{f}_{t,\ell} =\frac{1}{\lambda_{t}} ( w_{t,\ell}  \circ    \alpha_{t,\ell}^{-1}  )^{-1} \cdot  ( f_{t,\ell}  \circ    \alpha_{t,\ell}^{-1}  )^{-1}   .
				\]

				Set  
				\[
				A_{t,\ell} = \{ x \in K \;|\;  f_{t,\ell}(x) - 1 | \geq \varepsilon \}.
				\]
				
				By the same arguments as in the proof of \cite[Theorem 2.6]{arxiv2}    
				we can show that

				\[
				\parallel \chi_{A_{t,\ell}} \parallel_ \mathcal{Y}<  \epsilon 
				\qquad 
				\text{for all } t \in F,\; \ell \in \{1,\dots,N\}.
				\]
				
				Next, put
				\[
				C_{t,\ell}
				=
				\{ x \in K : 
				( \eta \circ \alpha_{t,\ell})(x) \cdot 
				w_{t,\ell}^{-1}(x) \cdot 
				\vert f_{t,\ell}(x) \vert
				>  m_K  \lambda_{t}  \varepsilon \}.
				\]
				
				Then, since
				\[
				\chi_{C_{t,\ell}}
				\bigl( \alpha_{t,\ell}^{-1}(K) \bigr) = 0
				\qquad 
				\text{for all } t \in F, \text{ and }  \ell \in \{1,\dots,N\}
				\]
				(because 
				$\alpha_{t,\ell}^{-1}(K) \cap K = \varnothing$),
				we get
				
				\[
				\frac{1}{\lambda_{t}} \chi_{C_{t,\ell}}
				( \eta \circ \alpha_{t,\ell} )\,
				w_{t,\ell}^{-1} \, \big| f_{t,\ell}
				\big|
				=
				\]
				
				\[
				=
				\bigl|
				\chi_{C_{t,\ell}}
				( \eta \circ \alpha_{t,\ell}) 
				\, [ ( \lambda_{t} w_{t,\ell} )^{-1} f_{t,\ell}
				- \chi_{K} \circ \alpha_{t,\ell}]
				\bigr|
				\le
				\]
				
				\[
				\le
				\bigl|
				(\eta \circ \alpha_{t,\ell})
				\bigl[
				( \lambda_{t} w_{t,\ell} )^{-1} f_{t,\ell}
				- \chi_{K} \circ \alpha_{t,\ell}
				\bigr]
				\bigr|
				=
				\]
				
				\[
				=
				\bigl|
				\Bigl[ 
				\eta \big(
				\frac{f_{t,\ell} \circ \alpha_{t,\ell}^{-1}}
				{\lambda_{t}\, ( w_{t,\ell} \circ \alpha_{t,\ell}^{-1} ) }
				- \chi_{K}\big)
				\Bigr]
				\circ \alpha_{t,\ell}
				\bigr|
				\]
				
				\[
				=
				\bigl|\eta \cdot
				(\widetilde{f}_{t,\ell} - \chi_{K})
				\bigr|\circ \alpha_{t,\ell}.
				\]
				
				Since
				
				\[
				\chi_{C_{t,\ell}} m_K \varepsilon
				<
				\frac{1}{\lambda_t}
				\chi_{C_{t,\ell}}
				( \eta_{t,\ell} \circ \alpha_{t,\ell} )
				\, w_{t,\ell}^{-1} \vert f_{t,\ell} \vert,
				\]
				
				by the same arguments as in the proof of \cite[Theorem 2.6]{arxiv2}   
				we can deduce that  
				
				\[
				\parallel 
				\chi_{C_{t,\ell}} \parallel  < \epsilon
				\]
				
				and  
				
				\[
				(\eta \circ \alpha_{t,\ell})(x)\,
				w_{t,\ell}^{-1}(x)\,
				<
				\frac{m_K \lambda_{t}  \varepsilon}
				{1 - \varepsilon}
				\]
				for all $t \in F,
				\ \ell \in \{1,\dots,N\}$
				and for all $x \in K \setminus (A_{t,\ell} \cup C_{t,\ell}).$\\
				Further, for each $t \in F$ and $\ell \in \{1,\dots,N\}$, put
				\[
				B_{t,\ell} 
				:=
				\{ x \in K : | \tilde{f}_{t}(x_{\ell}) - 1 | \ge \varepsilon \}.
				\]
				
				Then, by similar arguments as for $A_{t,\ell}$, we can show that
				\[
				\| \chi_{B_{t,\ell}} \|_{\mathcal{Y}} < \varepsilon
				\qquad
				\text{for all } t \in F,\; \ell \in \{1,\dots,N\}.
				\]
				Next, set
				\[
				D_{t,\ell}
				=
				\Bigl\{
				x \in K :
				(\eta \circ \alpha^{-1}_{t,\ell})(x)\,
				\cdot 
				(w^{-1}_{t,\ell} \circ  \alpha^{-1}_{t,\ell}) (x) \,
				|\tilde{f}_{t,\ell}(x)|
				>
				\frac{m_K \varepsilon}{\lambda_t}
				\Bigr\}.
				\]
				
				By similar arguments as for $C_{t,\ell}$, we can show that
				\[
				\lambda_t 
				\chi_{D_{t,\ell}}
				(\eta \circ \alpha^{-1}_{t,\ell}) (w_{t,\ell} \circ \alpha_{t,\ell}^{-1})
				\, |\tilde{f}_{t,\ell}|
				\le
				\Bigl|
				\Bigl(
				\eta (\lambda_t T_{t,\ell}(\tilde{f}_{t,\ell})
				- \chi_{K}
				\Bigr)
				\circ
				\alpha^{-1}_{t,\ell}
				\Bigr|,
				\]
				
				so that
				\[
				\| \chi_{D_{t,\ell}} \|_{\mathcal{Y}} < \varepsilon
				\qquad
				\text{and for all } x \in K \setminus (B_{t,\ell} \cup D_{t,\ell})
				\]
				
				we thus obtain
				\[
				(\eta \circ \alpha^{-1}_{t,\ell})(x)\,
				\cdot (w_{t,\ell}     \circ \alpha^{-1}_{t,\ell} )(x) 
				<
				\frac{m_K \varepsilon}{\lambda_t(1-\varepsilon)}.
				\]
				
				This holds for all $t \in F$ and $\ell \in \{1,\dots,N\}$.
				
				Since for each distinct $s,\ell \in \{1,\dots,N\}$ and all $t \in F$
				we have
				
				\[
				\lambda_t^{-1} T_{t,\ell}^{-1}\bigl( B(\chi_{K}, m_K \varepsilon^{2}) \bigr)
				\;\cap\;
				\lambda_t^{-1} T_{t,s}^{-1}\bigl( B(\chi_{K}, m_K \varepsilon^{2}) \bigr)
				\neq \varnothing,
				\]
				
				then
				
				\[
				T_{t,\ell}^{-1}\bigl( B(\chi_{K}, m_K \varepsilon^{2}) \bigr)
				\cap
				T_{t,s}^{-1}\bigl( B(\chi_{K}, m_K \varepsilon^{2}) \bigr)
				\neq \varnothing.
				\]
				
				So for each distinct $s,\ell \in \{1,\dots,N\}$ and all $t \in F$
				we can find some 
				\[
				h_{\ell,s,t} \in \mathcal{Y}_{\eta}
				\quad\text{such that}\quad
				\| T_{t,\ell} h_{\ell,s,t} - \chi_{K} \|_{\mathcal{Y}_{\eta}}
				< m_K \varepsilon^{2}.
				\]
				and $$ \| T_{t, s} h_{\ell,s,t} - \chi_{K} \|_{\mathcal{Y}_{\eta}}
				< m_K \varepsilon^{2}. $$ 
				
				Put
				\[
				\tilde{h}_{\ell,s,t}
				=T_{t, s} h_{\ell,s,t}
				\]
				Then $$
				h_{\ell,s,t} =\bigl( w_{t,s} \circ \alpha^{-1}_{t,s} \bigr)^{-1} \cdot
				\left(
				\tilde{h}_{\ell,s,t} \circ \alpha^{-1}_{t,s}
				\right).
				$$
				
				For each distinct $s,\ell \in \{1,\dots,N\}$ and $t \in F$, set
				
				\[
				H_{t,\ell,s}
				:=
				\Bigl\{
				x \in K :
				(\eta \circ \alpha^{-1}_{t,\ell}\circ \alpha_{t,s})(x)
				\,
				\frac{
					(w_{t,\ell} \circ \alpha^{-1}_{t,\ell}\circ \alpha_{t,s})(x)\,
					\vert	\tilde{h}_{\ell,s,t}(x) \vert
				}{
					w_{t,s}(x)
				}
				>
				m_K \varepsilon
				\Bigr\}.
				\]
				
				Then, since 
				\[
				\alpha_{t,\ell}^{-1} (\alpha_{t,s}(K)) \cap K = \varnothing 
				\qquad
				\text{for all } t \in F
				\text{ and distinct } s,\ell,
				\]
				we get
				\[
				\chi_{H_{t,\ell,s}}
				\bigl(
				\chi_{K} \circ \alpha^{-1}_{t,\ell} \circ \alpha_{t,s}
				\bigr)
				= 0.
				\]
				
				Hence,
				
				\[
				\Bigl|
				\chi_{H_{t,\ell,s}}
				(\eta \circ \alpha^{-1}_{t,\ell} \circ \alpha_{t,s})
				\frac{
					(w_{t,\ell} \circ \alpha^{-1}_{t,\ell} \circ \alpha_{t,s}\, )\, \tilde{h}_{\ell,s,t}
				}{
					w_{t,s}
				}
				\Bigr|
				\leq
				\]
				\[ \leq
				\Bigl|
				(\eta \circ \alpha^{-1}_{t,\ell}\circ \alpha_{t,s})
				\Bigl[
				\frac{
					(w_{t,\ell} \circ \alpha^{-1}_{t,\ell} \circ \alpha_{t,s}\, ) \tilde{h}_{\ell,s,t}
				}{
					w_{t,s}
				}
				-
				\chi_{K} \circ \alpha^{-1}_{t,\ell} \circ \alpha_{t,s}
				\Bigr]
				\Bigr|
				\]
				
				\[
				=
				\bigl|
				\eta \Bigl(
				T_{t,\ell}\bigl( (
				w_{t,s} \circ \alpha^{-1}_{t,s})
				^{-1} \cdot
				( \tilde{h}_{\ell,s,t}  \circ \alpha^{-1}_{t,s}) \bigr)
				-
				\chi_{K}
				\Bigr)
				\bigr| \circ \alpha^{-1}_{t,\ell}\circ \alpha_{t,s}
				\]
				
				\[
				=
				\bigl|
				\eta \bigl( T_{t,\ell} h_{\ell,s,t} - \chi_{K} \bigr)
				\bigr| \circ \alpha^{-1}_{t,\ell} \circ \alpha_{t,s}.
				\]
				
				By similar arguments as in the proof of \cite[Theorem 2.9]{arxiv2} ,  
				we can deduce that
				\[
				\|\chi_{H_{t,\ell,s}}\|_{\mathcal{Y}} < \varepsilon
				\qquad
				\text{for all } t \in F
				\text{ and each distinct } s,\ell \in \{1,\dots,N\}.
				\]
				
				Further, for each $t \in F$ and each distinct $s,\ell \in \{1,\dots,N\}$, set
				
				\[
				G_{t,s,\ell}
				:=
				\{ x \in K : | \tilde{h}_{\ell,s,t}(x) - 1 | \ge \varepsilon \}.
				\]
				
				Since
				\[
				\| \tilde{h}_{\ell,s,t} - \chi_{K} \|_{\mathcal{Y}_{\eta}} < m_K \varepsilon^{2},
				\]
				then by similar arguments as for $A_{t,\ell}$,
				we can conclude that
				\[
				\|\chi_{G_{t,s,\ell}}\|_{\mathcal{Y}} < \varepsilon.
				\]
				
				Moreover, for every $ x \in K \setminus (H_{t,\ell,s} \cup G_{t,\ell,s} ) $ we have that $$ (\eta \circ \alpha^{-1}_{t,\ell}\circ \alpha_{t,s})(x)
				\,
				\frac{
					(w_{t,\ell} \circ \alpha^{-1}_{t,\ell}\circ \alpha_{t,s})(x)\,
				}{
					w_{t,s}(x)
				} \leq \frac{m_K \varepsilon}{1-\varepsilon} .$$
				This holds for each $t \in F$ and every distinct $s,\ell \in \{1,\dots,N\} .$\\
				For each $t \in F$ put $$ E_t := K \setminus \Bigl( \bigcup_{1 \leq \ell, s \leq N, \ell \neq s} (A_{t,\ell} \cup B_{t,\ell} \cup C_{t,\ell} \cup D_{t,\ell} \cup H_{t,\ell,s} \cup G_{t,\ell,s} ) \Bigr) .$$ 
				
				Now we prove $(2)\Rightarrow(1)$.\\
				To this end, choose non-empty open subsets
				$\mathcal{O}, \mathcal{O}_{1}, \dots, \mathcal{O}_{N}
				\text{ of } \mathcal{Y}_{\eta}.$
				As in the proof of 
				$(2) \Rightarrow (1)$ in \cite[Theorem 2.6]{arxiv2} , we can find some $\widetilde{f}, \widetilde{g}_1, \dots, \widetilde{g}_N \in \mathcal{Y}_{bc}$ such that 
				$ \tilde{\tilde{f}} := \widetilde{f} \cdot \eta^{-1} \in \mathcal{O}$ and $\widetilde{\widetilde{g}}_\ell := \widetilde{g}_\ell \cdot  \eta^{-1} \in \mathcal{O}_\ell$ 
				for each $\ell \in \{1, \dots, N\}.$ There exists some positive constant $ \delta $ such that $ \delta-$neighbourhood of $ \tilde{\tilde{f}}$ is contained in $ \mathcal{O} $ and $ \delta-$neighbourhood of $ \widetilde{\widetilde{g}}_\ell$ is contained in $\mathcal{O}_\ell $ for each $\ell \in \{1, \dots, N\}.$
				
				Let
				\[
				K := \operatorname{supp}\tilde{f}\,
				\cup
				\Bigl( \bigcup_{\ell=1}^{N} \operatorname{supp}\,\tilde{g}_{\ell} \Bigr),
				\]
				and set $$ C := \max \lbrace \sup_{x \in X} (\widetilde{f}(x)),\sup_{x \in X} (\widetilde{g}_1 (x)), \cdots, \sup_{x \in X} (\widetilde{g}_N (x)) \rbrace  .$$ 
				Put $$ \varepsilon := \min \lbrace \frac{1}{2}, \frac{ \delta }{ (4 + 2N)N ( m_K + C)} \rbrace  .$$ 
				Choose $F \in \mathcal{F}$ and the family
				$\{ E_{t} \}_{t\in F}$
				of Borel subsets of $K$
				satisfying the assumptions of (2) with respect to $K$ and
				$\varepsilon$.
				
				Then, for all $t \in F$ and every distinct $s,\ell \in \{1,\dots,N\}$  
				we obtain  $$ \parallel T_{t, \ell} \big((
				w_{t,s}^{-1} \circ \alpha^{-1}_{t,s})((\widetilde{\widetilde{g}}_s \chi_{E_t}) \circ \alpha^{-1}_{t,s}) \big) \parallel_{\mathcal{Y}_{\eta}}= $$ 	 
				\[
				= \parallel
				( \eta \circ \alpha_{t,\ell}^{-1}  \circ \alpha_{t,s} )
				(
				w_{t,\ell} \circ \alpha_{t,\ell}^{-1}  \circ \alpha_{t,s} 
				)
				w_{t,s}^{-1} \widetilde{\widetilde{g}}_s \chi_{E_{t}} \parallel_{\mathcal{Y}} 
				\]
				\[
				\le
				\sup_{x \in E_{t}}
				\Bigl( 
				( \eta \circ \alpha_{t,\ell}^{-1}  \circ \alpha_{t,s} (x) )
				(
				w_{t,\ell} \circ \alpha_{t,\ell}^{-1}  \circ \alpha_{t,s} (x)
				) 	w_{t,s}^{-1}(x)
				\Bigr) 
				\parallel 
				\tilde{g_{s}}
				\parallel_{\mathcal{Y}} 
				\frac{1}{m_{k}}.
				\]   
				In addition, for all $t \in F$ and $s,\ell \in \{1,\dots,N\}$ we have
				\[
				\bigl\|
				T_{t,\ell}
				( \tilde{\tilde{f}} 
				\chi_{E_t} )
				\bigr\|_{\mathcal{Y}}
				\le
				\frac{1}{m_K}
				\sup_{x \in E_t}
				\bigl( (\eta \circ \alpha^{-1}_{t,\ell})(x) (w_{t,\ell} \circ \alpha^{-1}_{t,\ell})(x)  
				\big)
				\parallel \tilde{f} 
				\parallel_{\mathcal{Y}}
				\]
				and  
				\[
				\parallel (
				w_{t,s}^{-1} \circ \alpha^{-1}_{t,s})((\widetilde{\widetilde{g}}_s \chi_{E_t}) \circ \alpha^{-1}_{t,s})  \parallel_{\mathcal{Y}_{\eta}} \leq  
				\frac{1}{m_K}
				\sup_{x \in E_t}
				\bigl(
				(\eta \circ \alpha_{t,s})(x)
				(w_{t,s}^{-1} (x)
				\bigr)
				\parallel  \tilde{g}_{s} \parallel_{\mathcal{Y}}.
				\]
				
				This follows from the proof of \cite[Theorem 2.6]{arxiv2} .

				Finally, for all $t \in F$ and $s \in \{1,\dots,N\}$  
				we have
				
				\[
				T_{t,s}
				\Bigl(
				(w_{t,s}^{-1} \circ \alpha^{-1}_{t,s})
				\bigl(
				(\tilde{\tilde{g_{s}}} \chi_{E_t})
				\circ \alpha^{-1}_{t,s}
				\bigr)
				\Bigr)
				=\tilde{\tilde{g_{s}}} \chi_{E_t}.
				\]

				For each $t \in F$ set
				\[
				v_t
				=\tilde{\tilde{f}}
				\chi_{E_t}
				+
				\dfrac{\sqrt{\sum_{l=1}^{N}   \sup_{x \in E_t}
						\bigl( (\eta \circ \alpha^{-1}_{t,\ell})(x) (w_{t,\ell} \circ \alpha^{-1}_{t,\ell})(x)  }\bigl)}  {\sqrt{\sum_{l=1}^{N}  \sup_{x \in E_t}
						\bigl( (\eta \circ \alpha_{t,\ell})(x) (w_{t,\ell}^{-1} )(x) \bigl)}   } 
				\sum_{s=1}^{N}
				(w_{t,s} \circ \alpha^{-1}_{t,s})^{-1}   
				\bigl(
				(\tilde{\tilde{g_{s}}} \chi_{E_t})
				\circ \alpha^{-1}_{t,s}
				\bigr)
				\]
				
				By the above inequalities, the triangle inequality and the choice of $ \varepsilon $ it is not hard to deduce that $ v_t \in \mathcal{O} $ and $$ \dfrac {\sqrt{\sum_{l=1}^{N}  \sup_{x \in E_t}
						\bigl( (\eta \circ \alpha_{t,\ell})(x) (w_{t,\ell}^{-1} )(x) \bigl)}   }{\sqrt{\sum_{l=1}^{N}   \sup_{x \in E_t}
						\bigl( (\eta \circ \alpha^{-1}_{t,\ell})(x) (w_{t,\ell} \circ \alpha^{-1}_{t,\ell})(x)  }\bigl)} T_{t,s} (v_t ) \in \mathcal{O}_{s} $$  for all $t \in F$ and $s \in \{1,\dots,N\} .$

			\end{proof}
			\begin{example}\label{zavrsetak}
				Let $ N=1, X = \mathbb{R} \text{ and }  \eta, \alpha $ be as in Example \ref{primer-tezina}. Let $\mathcal{F}$ be the family of all infinite subsets of $\mathbb{N}.$ For each $ n \in \mathbb{N} ,$  let $ \alpha_n(t) = t-n $ for all $ t \in \mathbb{R} .$ Let $ w :=1 ,$ and for each $ n \in \mathbb{N} $ set $ w_n = \prod_{j=0}^{n - 1} (w\circ\alpha^j) .$ For any $ p,q$ with $ 1 \leq q < p < \infty ,$ let $ \mathcal{Y} = \tilde{\mathcal{M}}_{q}^{p}( \mathbb{R}) $ (for the definition and details of the construction of the space $ \tilde{\mathcal{M}}_{q}^{p}( \mathbb{R}) ,$ see \cite[Section 4]{saw}). Then the conditions of Theorem \ref{vekt} are satisfied.
			\end{example}
			
			\section{Open question for further research}
			
			At the end of this paper, we leave it as an open question for further research to find other examples of non-unital normed algebras and Banach bimodules satisfying the conditions from Section 3.
			
			\text{ }
			
			\textbf{Acknowledgement}: I am grateful to Professor S. M. Tabatabaie for suggesting linear dynamics on solid spaces as the topic for my research and for introducing to me the relevant literature.

	\bibliographystyle{amsplain}

\end{document}